\documentclass[reqno]{amsart}

\usepackage[utf8]{inputenc}
\usepackage[T1]{fontenc}
\usepackage[british]{babel,isodate}
\cleanlookdateon
\usepackage{lmodern}
\usepackage{enumerate}
\usepackage{tikz}
\usepackage{tikz-cd}
\usetikzlibrary{arrows}
\usepackage{cancel}
\usetikzlibrary{babel}
\usepackage{lscape}
\usepackage{url}
\usepackage{hyperref}
\usepackage{float}
\usepackage{graphicx, import}
\usepackage{mathtools}
\usepackage{dsfont}

\usepackage{xargs}                      
\usepackage{xcolor} 
\definecolor{OliveGreen}{rgb}{0,0.6,0}
 
\usepackage[colorinlistoftodos, prependcaption,textsize=tiny]{todonotes}
\newcommandx{\unsure}[2][1=]{\todo[linecolor=red,backgroundcolor=red!25,bordercolor=red,#1]{#2}}
\newcommandx{\change}[2][1=]{\todo[linecolor=blue,backgroundcolor=blue!25,bordercolor=blue,#1]{#2}}
\newcommandx{\info}[2][1=]{\todo[linecolor=OliveGreen,backgroundcolor=OliveGreen!25,bordercolor=OliveGreen,#1]{#2}}

\usepackage{amsmath,amsthm,amssymb, amsfonts, amscd}
\usepackage{rotating,subcaption}
\usepackage{url}
\usepackage{graphicx,bm}
\usepackage{mathtools}
\usepackage[mathscr]{euscript}
\allowdisplaybreaks
\usepackage{lipsum}

\usepackage{xpatch}
\makeatletter
\AtBeginDocument{\xpatchcmd{\@thm}{\thm@headpunct{.}}{\thm@headpunct{}}{}{}}
\makeatother

\newtheorem{theorem}{Theorem}[section]
\newtheorem{proposition}[theorem]{Proposition}
\newtheorem{lemma}[theorem]{Lemma}
\newtheorem{corollary}[theorem]{Corollary}

\theoremstyle{definition}
\newtheorem{definition}[theorem]{Definition}
\newtheorem{example}[theorem]{Example}

\newtheorem{construction}[theorem]{Construction}

\newtheorem{question}[theorem]{Question}

\theoremstyle{remark}

\newtheorem{remark}[theorem]{Remark}
\newtheorem{notation}[theorem]{Notation}

\numberwithin{equation}{section}

\newcommand{\hooklongrightarrow}{\lhook\joinrel\longrightarrow}
\newcommand{\N}{\mathbb{N}}

\newcommand{\uno}{\mathds{1}}

\newcommand{\C}{\mathcal{C}}

\newcommand{\id}{\mathrm{Id}}

\let\hom\relax
\newcommand{\hom}[3]{\mathrm{Hom}_{#1}(#2,#3)}
\renewcommand{\to}{\longrightarrow}

\newcommand{\toiso}{\overset{\cong}{\to}}

\begin{document}


\title{A combinatorial PROP for bialgebras}
\date{\today}

\author{Jorge Becerra}
\address{Bernouilli Institute, University of Groningen, The Netherlands}
\email{\href{mailto:j.becerra@rug.nl}{j.becerra@rug.nl}}
\urladdr{ \href{https://sites.google.com/view/becerra/}{https://sites.google.com/view/becerra/}} 




\begin{abstract}
It is a classical result that the category of finitely-generated free monoids serves as a PROP for commutative bialgebras. Attaching permutations to fix the order of multiplication, we construct an extension of this ca\-te\-gory that is equivalent to the PROP for bialgebras.
\end{abstract}


\maketitle

\setcounter{tocdepth}{1}
\tableofcontents


\section{Introduction}

Bialgebras are an important algebraic structure that one commonly finds in many areas of mathematics, eg algebraic topology \cite{maymoreconcise}, homotopical algebra \cite{lodayvallete}, quantum groups theory \cite{kassel}, knot theory \cite{habiro}, etc.  Therefore it is sometimes useful to take one level of abstraction up and study their structure maps on their own. A convenient setup to do this is given by symmetric monoidal categories whose family of objects can be identified with the set of non-negative integers, also known as \textit{PROPs}. One could think of a PROP as a gadget to study abstract algebraic operations, whereas an \textit{algebra} over a PROP is a concrete realisation of these operations. This realisation is achieved by means of a strong monoidal functor $\mathsf{P} \to \C$, where $\mathsf{P}$ is a PROP and $\C$ is a symmetric monoidal category.

In general, an algebraic structure is given by some structure maps that satisfy certain compatibility conditions. A rather straightforward way to describe a PROP for such an algebraic structure is given by setting the structure maps as generators and taking the compatibility conditions as relations. In this paper we focus on the PROP $\mathsf{B}$ for (non-commutative, non-cocommutative) bialgebras, monoidally generated by the maps $\mu: 2 \to 1$, $\eta: 0 \to 1$, $\Delta: 1 \to 2$ and $\varepsilon: 1 \to 0$, subject to the usual relations that define a bialgebra (see Figure \ref{fig bialgebra axioms}). Though it may seem this gives rise to complex combinations of these maps, it turns out that any morphism in $\mathsf{B}$ has a unique normal form.

\begin{theorem}[Theorem \ref{thm normal form}]\label{thm normal form intro}
Every morphism $f: n \to m$ in $\mathsf{B}$ factors in a unique way as $$f= (\mu^{[q_1]} \otimes \cdots \otimes \mu^{[q_m]}) \circ P_{\sigma} \circ (\Delta^{[p_1]} \otimes \cdots \otimes  \Delta^{[p_n]})  $$ for some (unique) integers $p_1, \ldots , p_n, q_1, \ldots , q_m,s \geq 0$ and a permutation ${\sigma \in \Sigma_s}$, where the map $\mu^{[k]}: k \to 1$ (resp. $\Delta^{[k]}: 1 \to k$) stands for the iterated mul\-ti\-pli\-ca\-tion (resp. comultiplication).
\begin{figure}[h]
\small
 \centering
\def\svgwidth{0.45\textwidth}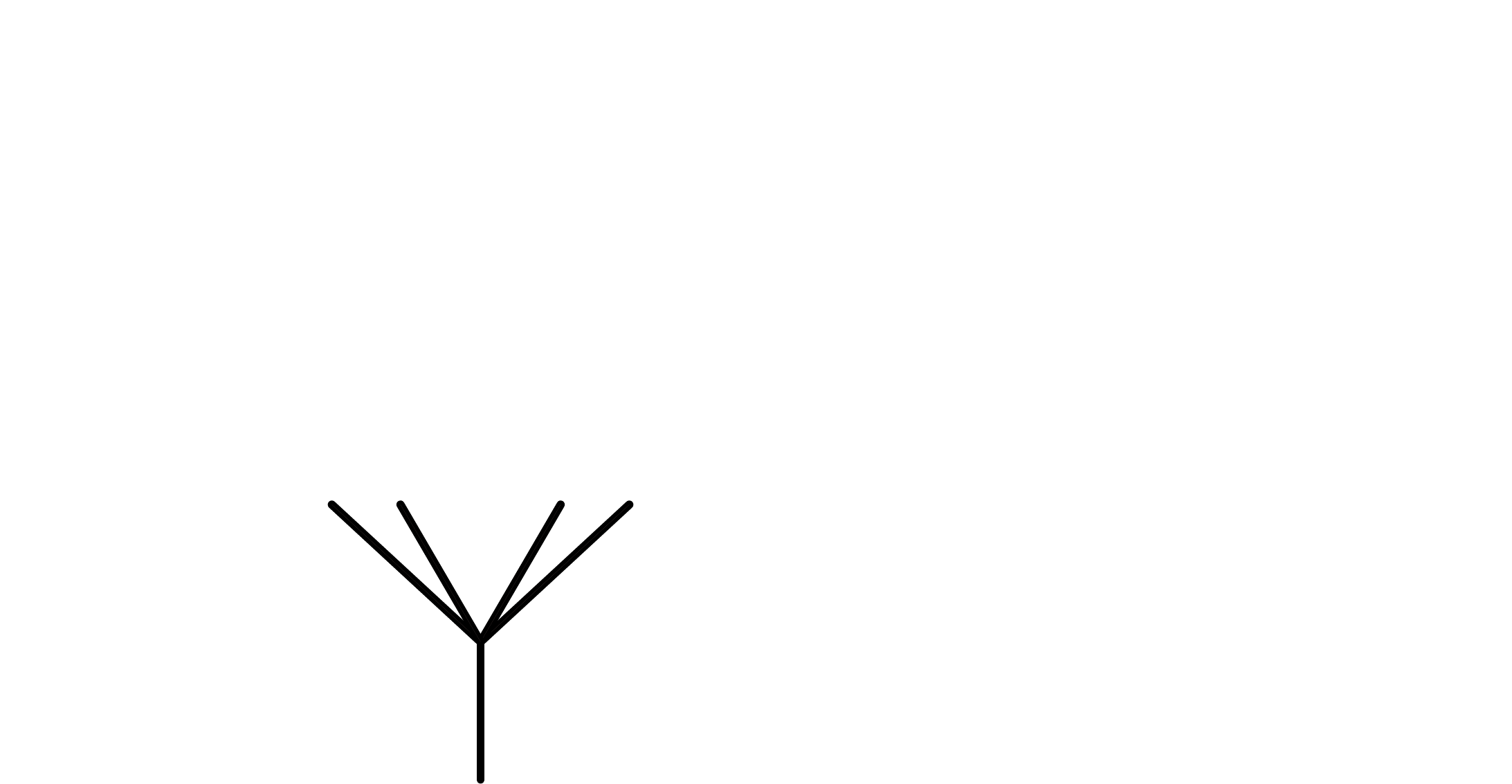
\end{figure}
\end{theorem}

Even though this is an important step towards the understanding of the PROP $\mathsf{B}$, computationally it does not present any advantage, as it gives no information about how to determine the coefficients $p_i, q_i, s$. Thence we would like to find an equivalent, computationally explicit description of $\mathsf{B}$.

Let  $\mathsf{fgFMon}$ be the category of finitely-generated free monoids and let $\C$ be a symmetric monoidal category. It turns out that the datum of a strong monoidal functor $\mathsf{fgFMon} \to \C$ uniquely determines a commutative bialgebra $A$ in $\C$, where the multiplication $\mu: A \otimes A \to A$ is given by the image of the monoid map $F(x,y) \to F(z)$, $x,y \mapsto z$ and the comultiplication $\Delta: A \to A \otimes A$ is the image of the monoid map $F(x) \to F(y,z)$, $x 	\mapsto yz$ (Example \ref{ex ComB = fgFMon}). In a similar fashion, if $\mathsf{fSet}$ denotes the category of finite sets, then the datum of a strong monoidal functor $\mathsf{fSet} \to \C$ uniquely determines a commutative algebra $A$ in $\C$, where this time the multiplication $\mu: A \otimes A \to A$ is given by the image of the unique map of sets $\{1,2 \} \to \{ 1 \}$ (Example \ref{ex fSet A}). These define equivalence of categories $\mathsf{ComA} \overset{\simeq}{\to} \mathsf{fSet} $  and $\mathsf{ComB} \overset{\simeq}{\to} \mathsf{fgFMon} $, where $\mathsf{ComA}$ and $\mathsf{ComB}$  are the PROPs for commutative algebras and bialgebras, respectively (we think of these as described in terms of generators and relations).


For instance, assume for simplicity $\C = \mathsf{Vect}_k$, the category of vector spaces over some field $k$. Given a commutative bialgebra $A$, consider the monoid maps $$f: F(x) \to F(a), \quad x \mapsto a^n \qquad , \qquad g:F(x) \to F(a,b), \quad x \mapsto aba^2b.$$ Under the previous equivalences, they  correspond to the $k$-linear maps
$$f': A \to A, \ x \mapsto x_{(1)} \cdots x_{(n)} \quad , \quad g': A \to A \otimes A, \ x  \mapsto x_{(1)} x_{(3)} x_{(4)} \otimes x_{(2)} x_{(5)}$$
where the order of the multiplication is arbitrary as $A$ is commutative. If $A$ was non-commutative, we would be able to fix the order of the multiplication by means of a permutation $\sigma \in \Sigma_n$  for $f'$ and a pair of permutations $(\tau_1, \tau_2) \in \Sigma_3 \times \Sigma_2$ for $g'$. The upshot of this observation is that by attaching these permutations to the monoid maps $f$ and $g$, we can produce non-commutative bialgebras from free monoids. A similar observation follows for finite set maps. 

One question that arises here is if this attaching of permutations can be made functorial in the following sense: there is a canonical faithful functor $\mathsf{ComA} \hooklongrightarrow \mathsf{ComB}$ with the property that precomposition with a commutative bialgebra object $\mathsf{ComB} \to \C$ yields its underlying commutative algebra object (for a precise definition of such a functor see \eqref{eq ComA Comb}). This embedding  corresponds, under the previous equivalences, to the free monoid functor, so that the diagram
$$\begin{tikzcd}
\mathsf{ComA} \dar{\simeq} \rar[hook] &  \mathsf{ComB}  \dar{\simeq}\\
\mathsf{fSet} \rar{F} & \mathsf{fgFMon}
\end{tikzcd}   $$
commutes (Lemma \ref{lem fSet fgFMon}). One would like to lift this diagram to a non-commutative context along the canonical full functors $\mathsf{A} \to \mathsf{ComA}$ and $\mathsf{B} \to \mathsf{ComB}$, where $\mathsf{A}$ is the PROP for (non-commutative) algebras.


Our main result states that the attaching of permutations gives an affirmative answer to the previous question.

\begin{theorem}[Theorem \ref{thm main thm 2}]\label{thm main thm}
There exist symmetric monoidal categories $\widehat{\mathsf{fSet}}$ and $\widehat{\mathsf{fgFMon}} $, together with full and essentially surjective functors $$\widehat{\mathsf{fSet}} \to \mathsf{fSet} \qquad , \qquad \widehat{\mathsf{fgFMon}} \to \mathsf{fgFMon},$$monoidal equivalences $$\mathsf{A} \overset{\simeq}{\to}\widehat{\mathsf{fSet}} \qquad , \qquad \mathsf{B} \overset{\simeq}{\to} \widehat{\mathsf{fgFMon}} $$
and a strong monoidal functor between them $$\widehat{F}: \widehat{\mathsf{fSet}} \to  \widehat{\mathsf{fgFMon}}  $$ making the following cube commutative,
\begin{equation}
\begin{tikzcd}[row sep=scriptsize,column sep=scriptsize]
& \widehat{\mathsf{fSet}} \arrow[from=dl,"\simeq"] \arrow[rr,"\widehat{F}"] \arrow[dd] & & \widehat{\mathsf{fgFMon}}  \arrow[from=dl,"\simeq"]\arrow[dd] \\
\mathsf{A} \arrow[rr,crossing over,hook]\arrow[dd] & & \mathsf{B} \\
& \mathsf{fSet}\arrow[from=dl, "\simeq"]\arrow[rr, "F \phantom{holaaa}"] &  & \mathsf{fgFMon}\arrow[from=dl,"\simeq"] \\
\mathsf{ComA} \arrow[rr,hook] & & \mathsf{ComB}\arrow[from=uu,crossing over]
\end{tikzcd}
\end{equation}
where $F$ is the free monoid functor.
\end{theorem}

The composite law in $\widehat{\mathsf{fgFMon}}$ gives an explicit formula that keeps track of the changes of indices  that arise when iterating the structure maps of a bialgebra, and that might be of interest in computational algebra, see \eqref{eq composite law fgFMonhat}. The category $\widehat{\mathsf{fSet}}$ is a suitable modification of the category $\widetilde{\mathsf{fSet}}$ of finite sets and maps of sets with a total order on fibres (Proposition \ref{prop fSethat fSettilde}) originally described by Pirashvili \cite{pirashvili} under the name of $\mathcal{F}(\mathrm{as})$. On the other hand, the category $\widehat{\mathsf{fgFMon}}$ can be viewed as a  combinatorial, explicit description of the Quillen $Q$-construction $\mathcal{Q}\mathcal{F}(\mathrm{as})$ carried out by Pirashvili. In particular, the composite law of $\widehat{\mathsf{fgFMon}}$ encodes an expression for the bialgebra maps $1 \to 1$ from \cite[5.3]{pirashvili} as a special case.

\subsection*{Outline of the paper} The paper is structured as follows: in Section \ref{sec Bialgebras and PROPs} we recall the definition of a bialgebra and fix notation. We also recall the notion of PROP and algebra over a PROP and a general construction via generators and relations. Examples \ref{ex fSet A} and \ref{ex ComB = fgFMon} make explicit the equivalences $\mathsf{ComA} \simeq \mathsf{fSet}$ and  $\mathsf{ComB} \simeq \mathsf{fgFMon}$, and in Proposition \ref{prop fSethat fSettilde} we lift the previous equivalence for algebras to the non-commutative setting, showing that $\widehat{\mathsf{fSet}} \simeq \mathsf{A}$.

In Section \ref{sec A PROP for bialgebras} we start by a few motivating examples showing the necessity to choose some permutations to fix the order of multiplication.  Next we briefly detour to prove Theorem \ref{thm normal form intro}, which will be needed to construct a PROP equivalent to $\mathsf{B}$. We later define the category $\widehat{\mathsf{fgFMon}}$ and show the equivalence $\widehat{\mathsf{fgFMon}} \simeq \mathsf{B}$ in Theorem \ref{thm fundamental theorem}. The key step is to define properly the composite law in $\widehat{\mathsf{fgFMon}}$, for which we will need a few constructions relating permutations and elements of a free monoid. These are carried out in pages \pageref{constr Phi} -- \pageref{lem utilisimo}. We end the section with a couple of examples illustrating our construction. 

Lastly, in Section \ref{sec Applications} we include several rather immediate applications of the equivalence $\widehat{\mathsf{fgFMon}} \simeq \mathsf{B}$ and prove Theorem \ref{thm main thm}.

\subsection*{Acknowledgments} The author would like to thank the anonymous referees for their extensive, detailed and thorough comments and suggestions, and for pointing out an error in a proof of the previous version of the manuscript. The author would also like to thank for Roland van der Veen for helpful discussions and valuable comments during the preparation of this paper.

\section{Bialgebras and PROPs}\label{sec Bialgebras and PROPs}

For this section we fix $(\C, \otimes, \mathds{1},P)$ a symmetric monoidal category, where $\otimes$ stands for the monoidal product, $\mathds{1}$ stands for the unit object and $P: \otimes \overset{\cong}{\Longrightarrow} \otimes^{op}$ stands for the symmetry, so $P_{Y, X} \circ P_{X,Y} = \id_{X\otimes Y}$. We will implicitly use MacLane's   coherence theorem \cite{maclane} to remove the associativity and unit constrains from the formulae. For an introduction to monoidal categories see \cite{turaevvirelizier}.

\subsection{Algebras, coalgebras and bialgebras}\label{subsec Algebras, coalgebras and bialgebras}
An \textit{algebra} in $\C$ (more precisely, an \textit{algebra object} in $\C$) is an object $A \in \C$ together with two arrows $\mu: A \otimes A \to A$ and $\eta: \uno \to A$, called the \textit{multiplication} and the \textit{unit}, satisfying the associativity and unit conditions
\begin{equation}\label{eq algebra axioms}
\mu \circ (\mu \otimes \id) = \mu \circ (\id \otimes \mu) \qquad , \qquad \mu \circ (\id \otimes \eta) = \id = \mu \circ (\eta \otimes \id).
\end{equation}
If $\mu \circ P_{A,A} =\mu$, we say that $A$ is \textit{commutative}.

A \textit{coalgebra} in $\C$ is an algebra in $\C^{op}$, that is, an object $A \in \C$ endowed with two arrows $\Delta:  A\to A \otimes A $ and $\varepsilon: A \to \uno $, called the \textit{comultiplication} and the \textit{counit}, satisfying the coassociativity and counit conditions
\begin{equation}\label{eq coalgebra axioms}
(\Delta \otimes \id)\circ \Delta =(\id \otimes \Delta) \circ \Delta \qquad , \qquad (\id \otimes \varepsilon)\circ \Delta = \id = (\varepsilon \otimes \id)\circ \Delta.  
\end{equation} 
If $ P_{A,A} \circ \Delta =\Delta$, we say that $A$ is \textit{cocommutative}.
It should be clear that a \textit{(co)algebra morphism} between (co)algebras should be an arrow in $\C$ that respects the structure maps.

If $A, A'$ are algebras in $\C$, the monoidal product $A \otimes A'$ is naturally an algebra with structure maps $(\mu \otimes \mu')(\id \otimes P_{A,A'} \otimes \id)$ and $\eta \otimes \eta'$. The same observation with reversed arrows follows for coalgebras. A \textit{bialgebra} in $\C$ is an  object $A \in \C$ which is both an algebra and a coalgebra and whose structure maps are compatible in the sense that the coalgebra structure maps are algebra morphisms (alternatively, the algebra structure maps are coalgebra morphisms). Explicitly, 

\begin{equation}\label{eq bialgebra axioms}
\begin{gathered}
\Delta \circ  \mu = (\mu \otimes \mu)\circ(\id \otimes P_{A,A} \otimes \id)\circ (\Delta \otimes \Delta) \\
\Delta \circ \eta = \eta \otimes \eta \quad , \quad \varepsilon \circ \mu = \varepsilon \otimes \varepsilon \quad , \quad \varepsilon \circ  \eta = \id_{\mathds{1}}.
\end{gathered}
\end{equation}

\begin{figure}[b]
\centering
\begin{subfigure}{.25\linewidth}
    \centering
\def\svgwidth{\textwidth}
\begingroup%
  \makeatletter%
  \providecommand\color[2][]{%
    \errmessage{(Inkscape) Color is used for the text in Inkscape, but the package 'color.sty' is not loaded}%
    \renewcommand\color[2][]{}%
  }%
  \providecommand\transparent[1]{%
    \errmessage{(Inkscape) Transparency is used (non-zero) for the text in Inkscape, but the package 'transparent.sty' is not loaded}%
    \renewcommand\transparent[1]{}%
  }%
  \providecommand\rotatebox[2]{#2}%
  \newcommand*\fsize{\dimexpr\f@size pt\relax}%
  \newcommand*\lineheight[1]{\fontsize{\fsize}{#1\fsize}\selectfont}%
  \ifx\svgwidth\undefined%
    \setlength{\unitlength}{682.79168557bp}%
    \ifx\svgscale\undefined%
      \relax%
    \else%
      \setlength{\unitlength}{\unitlength * \real{\svgscale}}%
    \fi%
  \else%
    \setlength{\unitlength}{\svgwidth}%
  \fi%
  \global\let\svgwidth\undefined%
  \global\let\svgscale\undefined%
  \makeatother%
  \begin{picture}(1,0.66633986)%
    \lineheight{1}%
    \setlength\tabcolsep{0pt}%
    \put(0,0){\includegraphics[width=\unitlength,page=1]{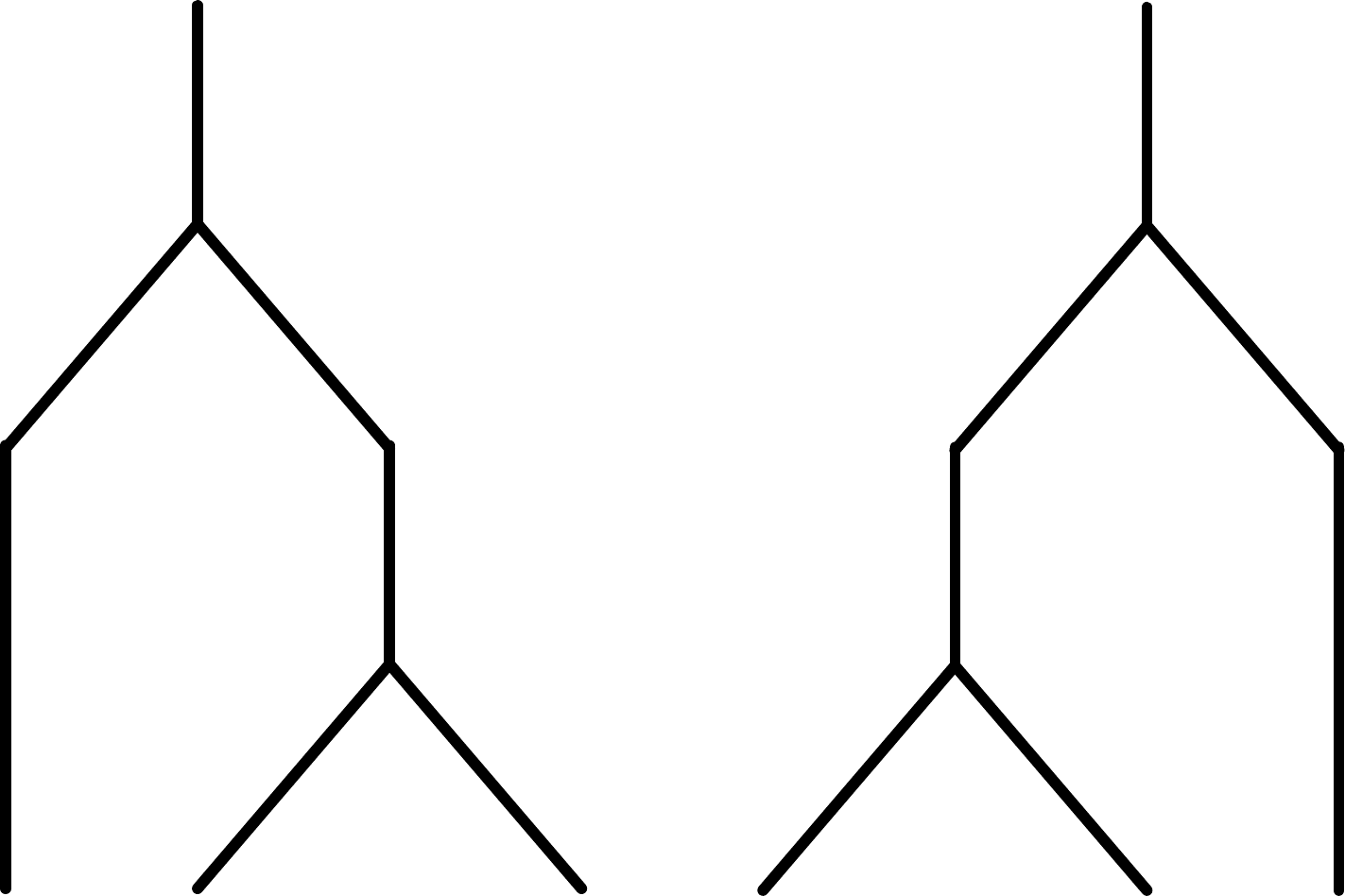}}%
    \put(0.44329477,0.29606083){\makebox(0,0)[lt]{\lineheight{1.25}\smash{\begin{tabular}[t]{l}$=$\end{tabular}}}}%
  \end{picture}%
\endgroup%

    \caption{}\label{fig associativity}
\end{subfigure}
    \hfill
\begin{subfigure}{.25\linewidth}
    \centering
\def\svgwidth{\textwidth}
\begingroup%
  \makeatletter%
  \providecommand\color[2][]{%
    \errmessage{(Inkscape) Color is used for the text in Inkscape, but the package 'color.sty' is not loaded}%
    \renewcommand\color[2][]{}%
  }%
  \providecommand\transparent[1]{%
    \errmessage{(Inkscape) Transparency is used (non-zero) for the text in Inkscape, but the package 'transparent.sty' is not loaded}%
    \renewcommand\transparent[1]{}%
  }%
  \providecommand\rotatebox[2]{#2}%
  \newcommand*\fsize{\dimexpr\f@size pt\relax}%
  \newcommand*\lineheight[1]{\fontsize{\fsize}{#1\fsize}\selectfont}%
  \ifx\svgwidth\undefined%
    \setlength{\unitlength}{767.91293407bp}%
    \ifx\svgscale\undefined%
      \relax%
    \else%
      \setlength{\unitlength}{\unitlength * \real{\svgscale}}%
    \fi%
  \else%
    \setlength{\unitlength}{\svgwidth}%
  \fi%
  \global\let\svgwidth\undefined%
  \global\let\svgscale\undefined%
  \makeatother%
  \begin{picture}(1,0.5965694)%
    \lineheight{1}%
    \setlength\tabcolsep{0pt}%
    \put(0,0){\includegraphics[width=\unitlength,page=1]{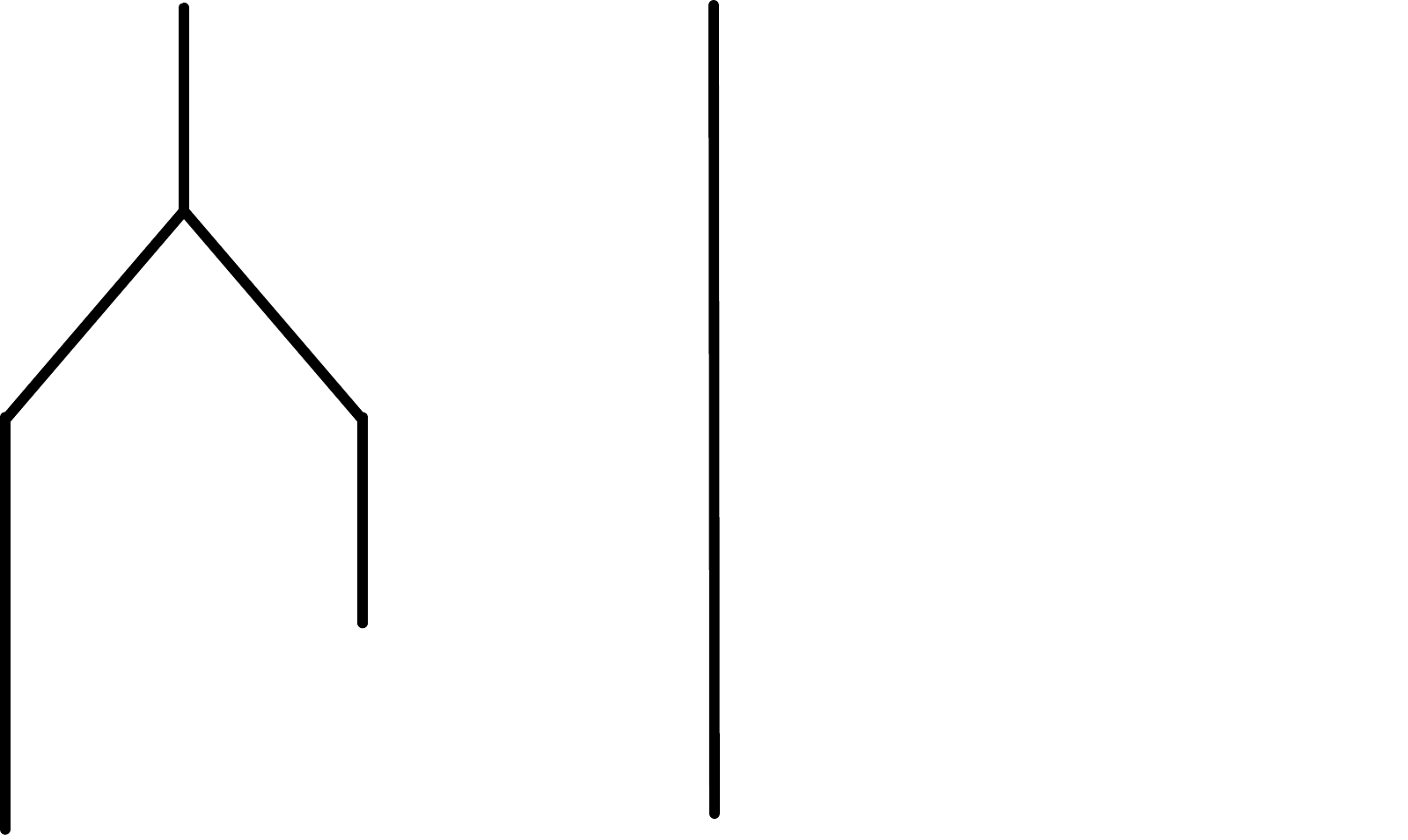}}%
    \put(0.34704298,0.29433915){\makebox(0,0)[lt]{\lineheight{1.25}\smash{\begin{tabular}[t]{l}$=$\end{tabular}}}}%
    \put(0,0){\includegraphics[width=\unitlength,page=2]{unitality_product.pdf}}%
    \put(0.56942453,0.29796951){\makebox(0,0)[lt]{\lineheight{1.25}\smash{\begin{tabular}[t]{l}$=$\end{tabular}}}}%
  \end{picture}%
\endgroup%

    \caption{}\label{fig unitality}
\end{subfigure}
   \hfill
\begin{subfigure}{.25\linewidth}
    \centering
\def\svgwidth{\textwidth}
\begingroup%
  \makeatletter%
  \providecommand\color[2][]{%
    \errmessage{(Inkscape) Color is used for the text in Inkscape, but the package 'color.sty' is not loaded}%
    \renewcommand\color[2][]{}%
  }%
  \providecommand\transparent[1]{%
    \errmessage{(Inkscape) Transparency is used (non-zero) for the text in Inkscape, but the package 'transparent.sty' is not loaded}%
    \renewcommand\transparent[1]{}%
  }%
  \providecommand\rotatebox[2]{#2}%
  \newcommand*\fsize{\dimexpr\f@size pt\relax}%
  \newcommand*\lineheight[1]{\fontsize{\fsize}{#1\fsize}\selectfont}%
  \ifx\svgwidth\undefined%
    \setlength{\unitlength}{717.75982714bp}%
    \ifx\svgscale\undefined%
      \relax%
    \else%
      \setlength{\unitlength}{\unitlength * \real{\svgscale}}%
    \fi%
  \else%
    \setlength{\unitlength}{\svgwidth}%
  \fi%
  \global\let\svgwidth\undefined%
  \global\let\svgscale\undefined%
  \makeatother%
  \begin{picture}(1,0.63373641)%
    \lineheight{1}%
    \setlength\tabcolsep{0pt}%
    \put(0,0){\includegraphics[width=\unitlength,page=1]{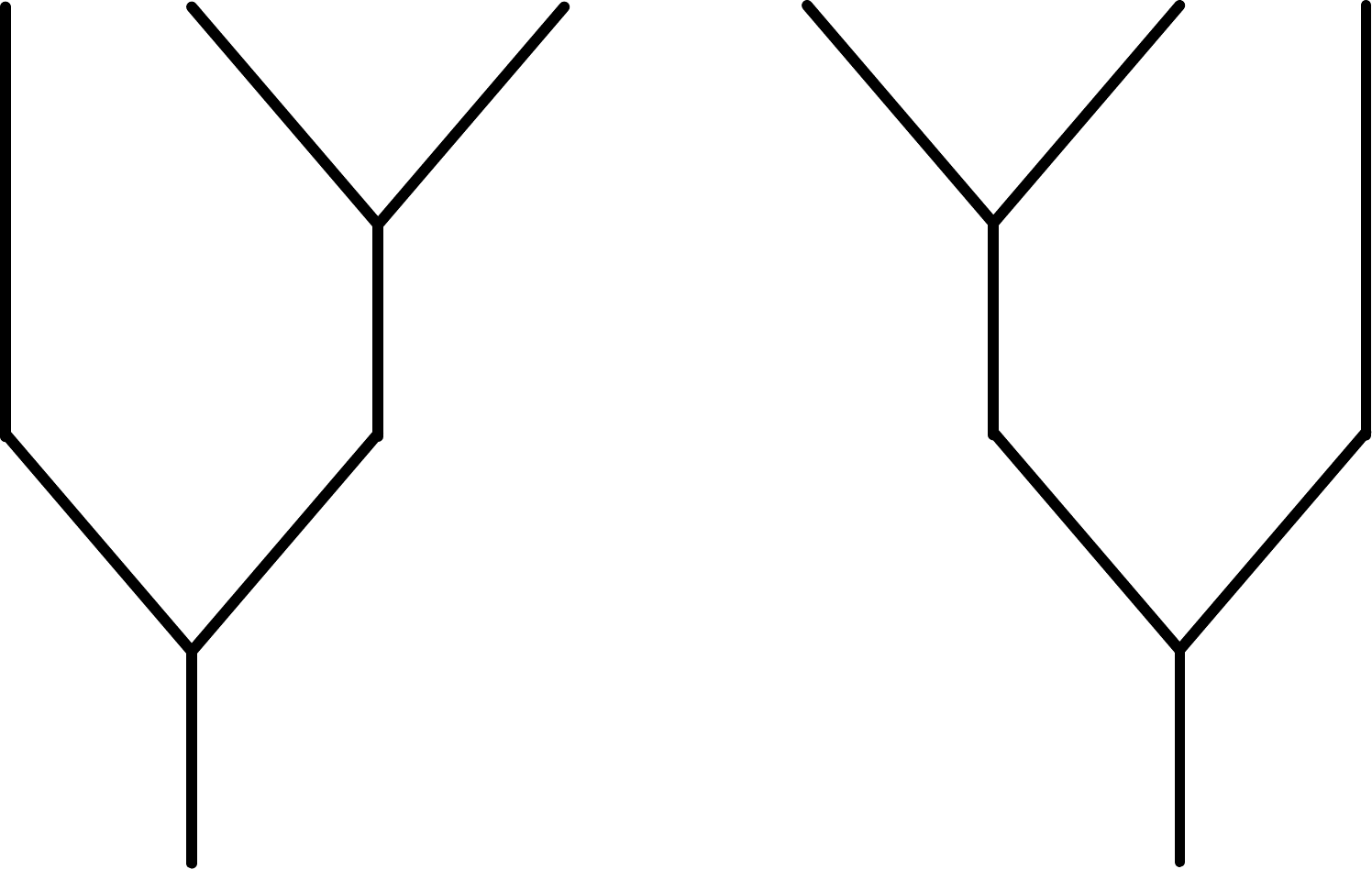}}%
    \put(0.46914315,0.32390857){\makebox(0,0)[lt]{\lineheight{1.25}\smash{\begin{tabular}[t]{l}$=$\end{tabular}}}}%
  \end{picture}%
\endgroup%

    \caption{}\label{fig coassociativity}
\end{subfigure}

\bigskip
\begin{subfigure}{0.25\linewidth}
  \centering
\def\svgwidth{\textwidth}
\begingroup%
  \makeatletter%
  \providecommand\color[2][]{%
    \errmessage{(Inkscape) Color is used for the text in Inkscape, but the package 'color.sty' is not loaded}%
    \renewcommand\color[2][]{}%
  }%
  \providecommand\transparent[1]{%
    \errmessage{(Inkscape) Transparency is used (non-zero) for the text in Inkscape, but the package 'transparent.sty' is not loaded}%
    \renewcommand\transparent[1]{}%
  }%
  \providecommand\rotatebox[2]{#2}%
  \newcommand*\fsize{\dimexpr\f@size pt\relax}%
  \newcommand*\lineheight[1]{\fontsize{\fsize}{#1\fsize}\selectfont}%
  \ifx\svgwidth\undefined%
    \setlength{\unitlength}{767.91293407bp}%
    \ifx\svgscale\undefined%
      \relax%
    \else%
      \setlength{\unitlength}{\unitlength * \real{\svgscale}}%
    \fi%
  \else%
    \setlength{\unitlength}{\svgwidth}%
  \fi%
  \global\let\svgwidth\undefined%
  \global\let\svgscale\undefined%
  \makeatother%
  \begin{picture}(1,0.5965694)%
    \lineheight{1}%
    \setlength\tabcolsep{0pt}%
    \put(0,0){\includegraphics[width=\unitlength,page=1]{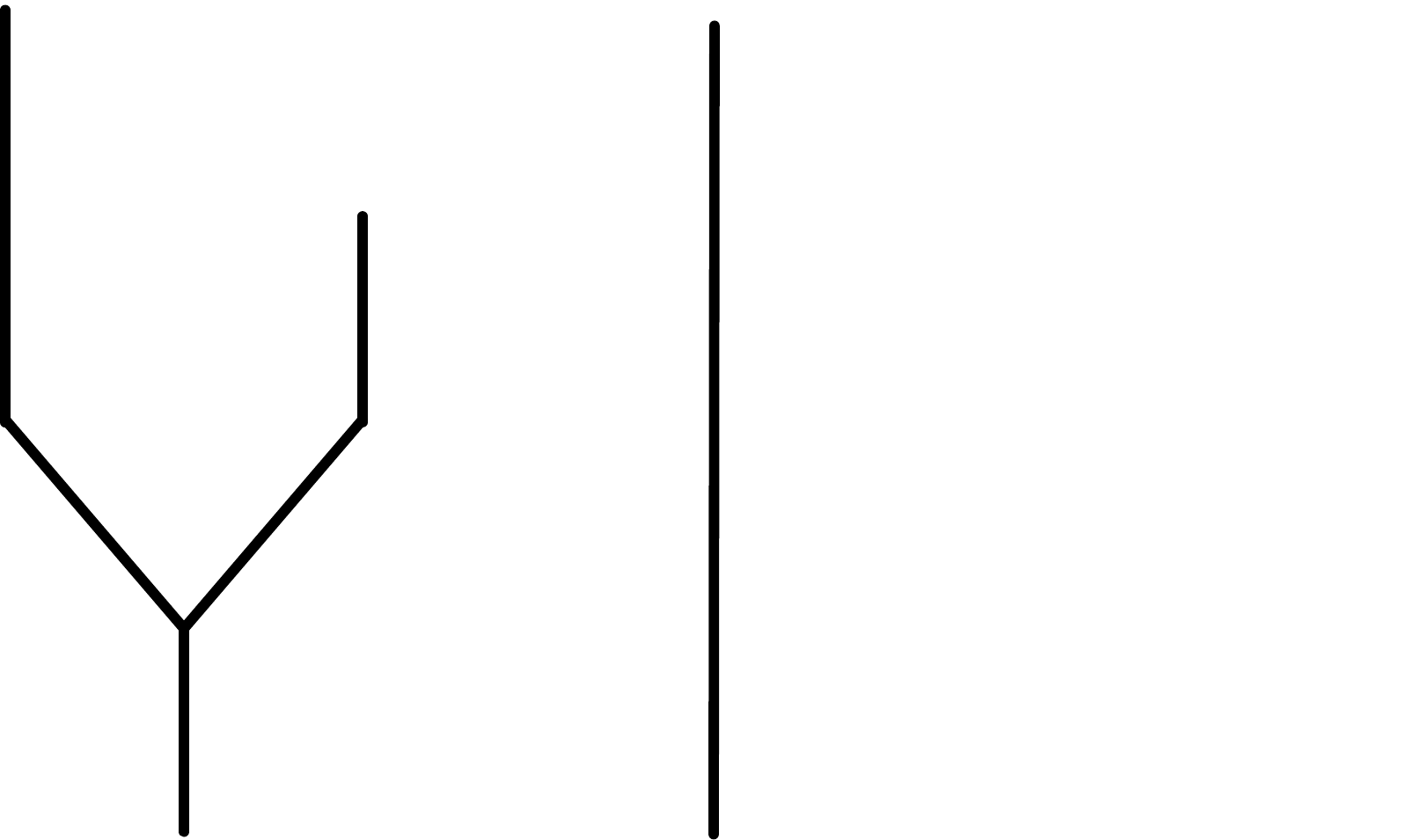}}%
    \put(0.34704298,0.27829037){\makebox(0,0)[lt]{\lineheight{1.25}\smash{\begin{tabular}[t]{l}$=$\end{tabular}}}}%
    \put(0,0){\includegraphics[width=\unitlength,page=2]{counitality_product.pdf}}%
    \put(0.56942453,0.27466){\makebox(0,0)[lt]{\lineheight{1.25}\smash{\begin{tabular}[t]{l}$=$\end{tabular}}}}%
  \end{picture}%
\endgroup%

  \caption{}\label{fig counitality}
\end{subfigure} 
   \hfill
\begin{subfigure}{0.25\linewidth}
  \centering
\def\svgwidth{\textwidth}
\begingroup%
  \makeatletter%
  \providecommand\color[2][]{%
    \errmessage{(Inkscape) Color is used for the text in Inkscape, but the package 'color.sty' is not loaded}%
    \renewcommand\color[2][]{}%
  }%
  \providecommand\transparent[1]{%
    \errmessage{(Inkscape) Transparency is used (non-zero) for the text in Inkscape, but the package 'transparent.sty' is not loaded}%
    \renewcommand\transparent[1]{}%
  }%
  \providecommand\rotatebox[2]{#2}%
  \newcommand*\fsize{\dimexpr\f@size pt\relax}%
  \newcommand*\lineheight[1]{\fontsize{\fsize}{#1\fsize}\selectfont}%
  \ifx\svgwidth\undefined%
    \setlength{\unitlength}{699.11499023bp}%
    \ifx\svgscale\undefined%
      \relax%
    \else%
      \setlength{\unitlength}{\unitlength * \real{\svgscale}}%
    \fi%
  \else%
    \setlength{\unitlength}{\svgwidth}%
  \fi%
  \global\let\svgwidth\undefined%
  \global\let\svgscale\undefined%
  \makeatother%
  \begin{picture}(1,0.65037947)%
    \lineheight{1}%
    \setlength\tabcolsep{0pt}%
    \put(0,0){\includegraphics[width=\unitlength,page=1]{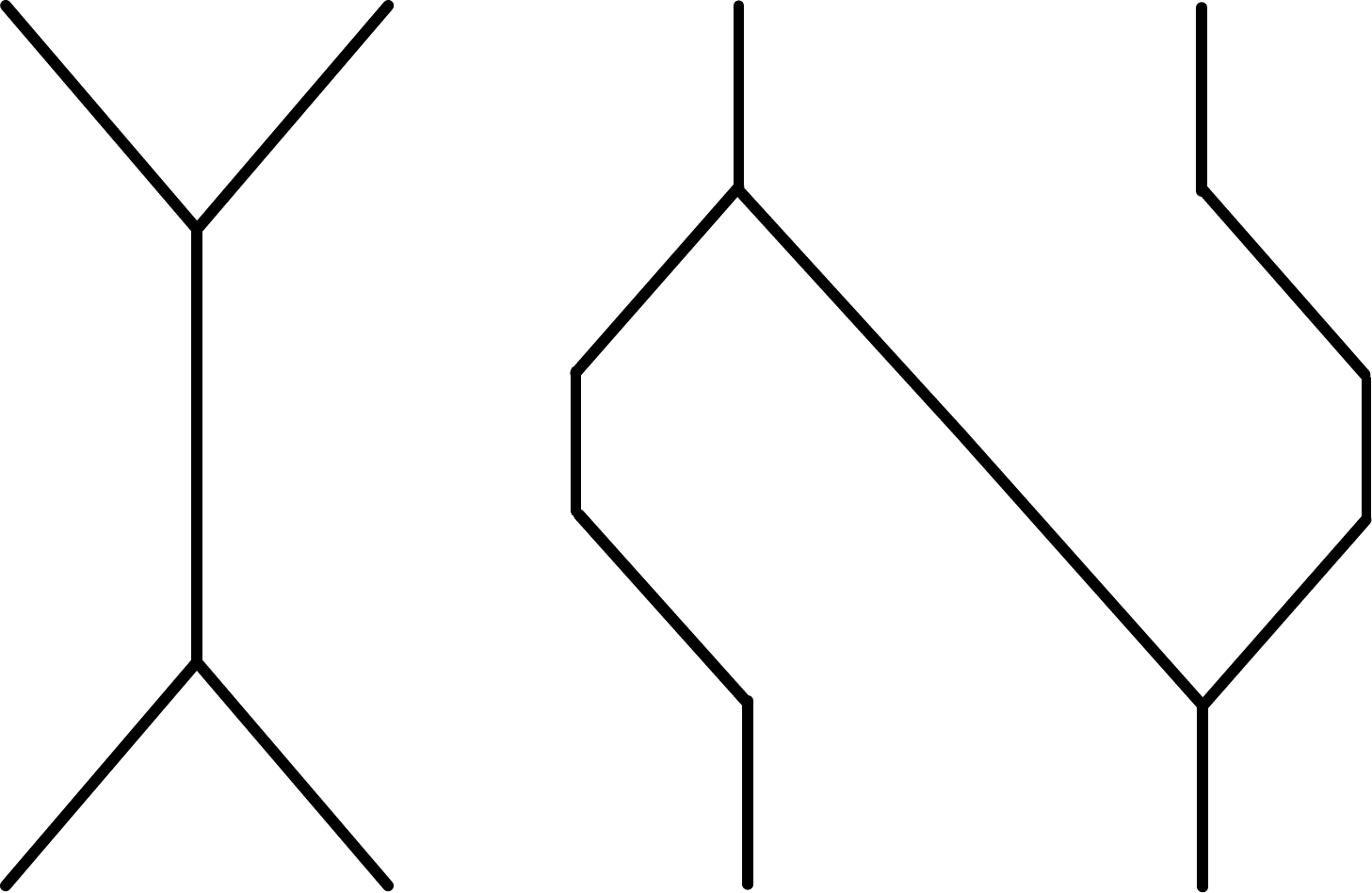}}%
    \put(0.2450883,0.32188103){\makebox(0,0)[lt]{\lineheight{1.25}\smash{\begin{tabular}[t]{l}$=$ \end{tabular}}}}%
    \put(0,0){\includegraphics[width=\unitlength,page=2]{bialgebra1.pdf}}%
  \end{picture}%
\endgroup%

  \caption{}\label{fig bialgebra1}
\end{subfigure} 
   \hfill
\begin{subfigure}{0.25\linewidth}
  \centering
\def\svgwidth{\textwidth}
\begingroup%
  \makeatletter%
  \providecommand\color[2][]{%
    \errmessage{(Inkscape) Color is used for the text in Inkscape, but the package 'color.sty' is not loaded}%
    \renewcommand\color[2][]{}%
  }%
  \providecommand\transparent[1]{%
    \errmessage{(Inkscape) Transparency is used (non-zero) for the text in Inkscape, but the package 'transparent.sty' is not loaded}%
    \renewcommand\transparent[1]{}%
  }%
  \providecommand\rotatebox[2]{#2}%
  \newcommand*\fsize{\dimexpr\f@size pt\relax}%
  \newcommand*\lineheight[1]{\fontsize{\fsize}{#1\fsize}\selectfont}%
  \ifx\svgwidth\undefined%
    \setlength{\unitlength}{699.11486047bp}%
    \ifx\svgscale\undefined%
      \relax%
    \else%
      \setlength{\unitlength}{\unitlength * \real{\svgscale}}%
    \fi%
  \else%
    \setlength{\unitlength}{\svgwidth}%
  \fi%
  \global\let\svgwidth\undefined%
  \global\let\svgscale\undefined%
  \makeatother%
  \begin{picture}(1,0.65037959)%
    \lineheight{1}%
    \setlength\tabcolsep{0pt}%
    \put(0,0){\includegraphics[width=\unitlength,page=1]{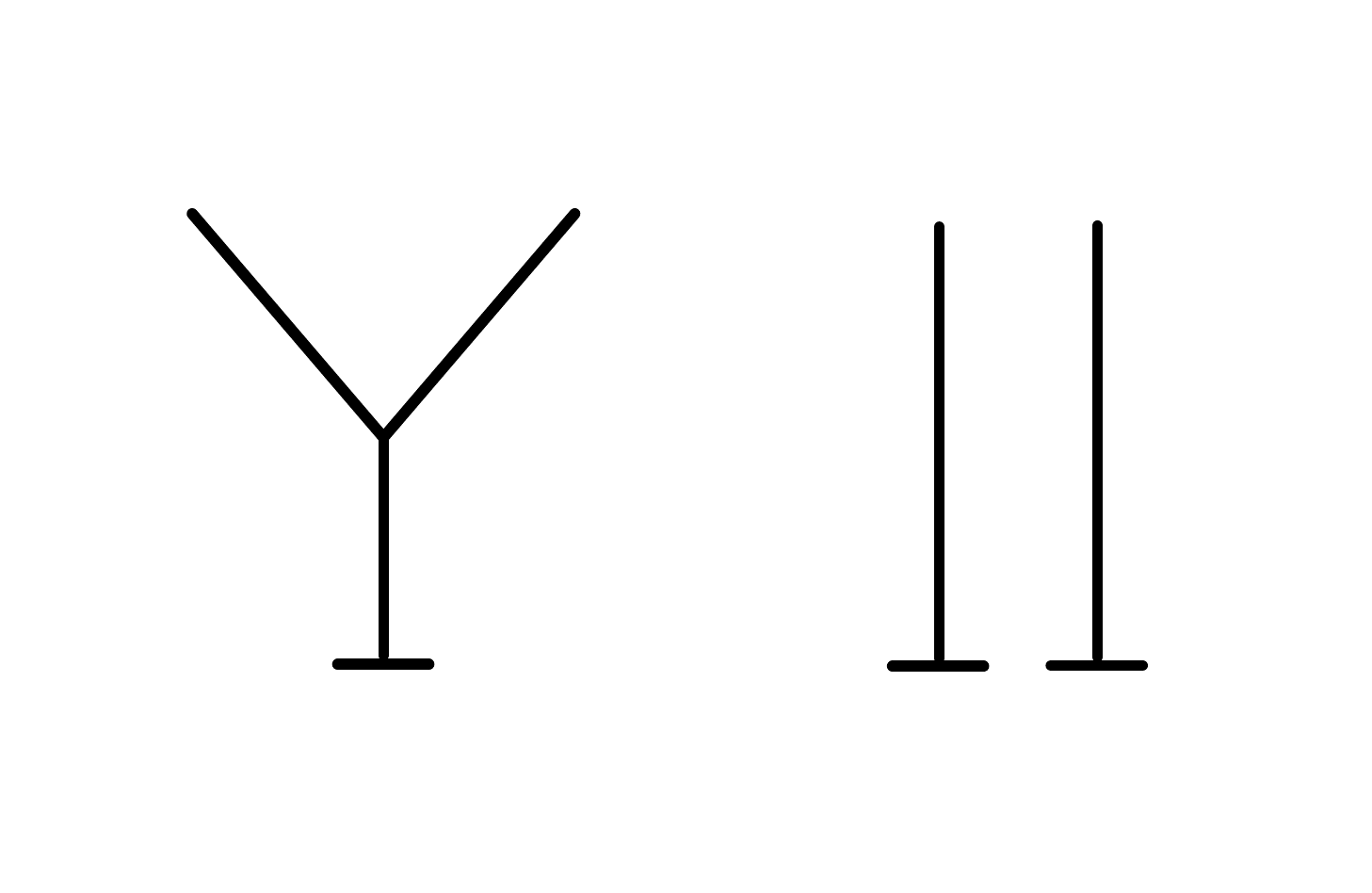}}%
    \put(0.46524565,0.3124238){\color[rgb]{0,0,0}\makebox(0,0)[lt]{\lineheight{1.25}\smash{\begin{tabular}[t]{l}$=$\end{tabular}}}}%
    \put(0,0){\includegraphics[width=\unitlength,page=2]{bialgebra2.pdf}}%
  \end{picture}%
\endgroup%

  \caption{}\label{fig bialgebra2}
\end{subfigure}

\bigskip
\begin{subfigure}{0.12\linewidth}
\phantom{.}
\end{subfigure} 
   \hfill
\begin{subfigure}{0.25\linewidth}
  \centering
\def\svgwidth{\textwidth}
\begingroup%
  \makeatletter%
  \providecommand\color[2][]{%
    \errmessage{(Inkscape) Color is used for the text in Inkscape, but the package 'color.sty' is not loaded}%
    \renewcommand\color[2][]{}%
  }%
  \providecommand\transparent[1]{%
    \errmessage{(Inkscape) Transparency is used (non-zero) for the text in Inkscape, but the package 'transparent.sty' is not loaded}%
    \renewcommand\transparent[1]{}%
  }%
  \providecommand\rotatebox[2]{#2}%
  \newcommand*\fsize{\dimexpr\f@size pt\relax}%
  \newcommand*\lineheight[1]{\fontsize{\fsize}{#1\fsize}\selectfont}%
  \ifx\svgwidth\undefined%
    \setlength{\unitlength}{699.11486047bp}%
    \ifx\svgscale\undefined%
      \relax%
    \else%
      \setlength{\unitlength}{\unitlength * \real{\svgscale}}%
    \fi%
  \else%
    \setlength{\unitlength}{\svgwidth}%
  \fi%
  \global\let\svgwidth\undefined%
  \global\let\svgscale\undefined%
  \makeatother%
  \begin{picture}(1,0.65037959)%
    \lineheight{1}%
    \setlength\tabcolsep{0pt}%
    \put(0,0){\includegraphics[width=\unitlength,page=1]{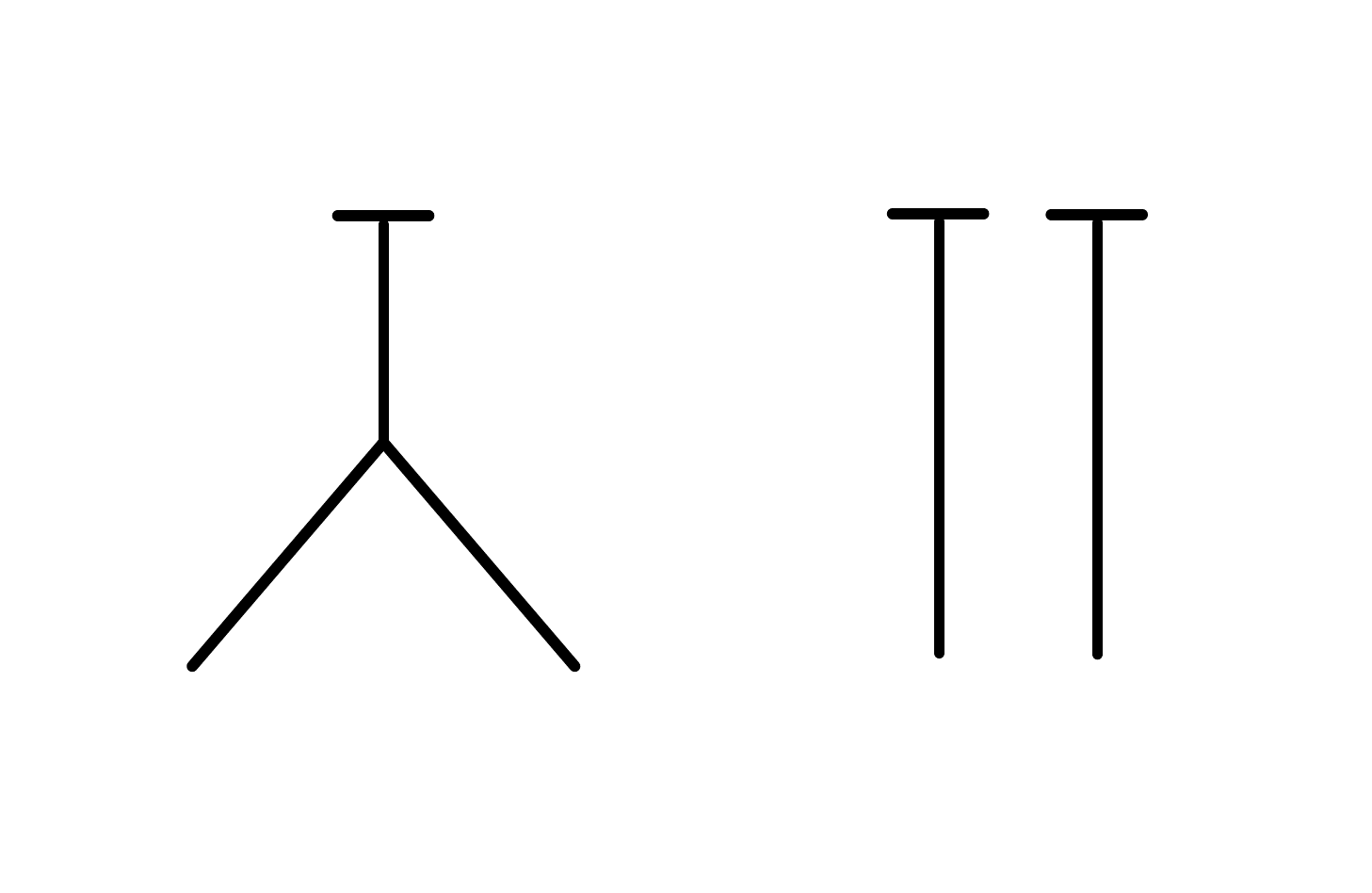}}%
    \put(0.46524565,0.32146057){\color[rgb]{0,0,0}\makebox(0,0)[lt]{\lineheight{1.25}\smash{\begin{tabular}[t]{l}$=$\end{tabular}}}}%
    \put(0,0){\includegraphics[width=\unitlength,page=2]{bialgebra3.pdf}}%
  \end{picture}%
\endgroup%

  \caption{}\label{fig bialgebra3}
\end{subfigure} 
   \hfill
\begin{subfigure}{0.25\linewidth}
  \centering
\def\svgwidth{\textwidth}
\begingroup%
  \makeatletter%
  \providecommand\color[2][]{%
    \errmessage{(Inkscape) Color is used for the text in Inkscape, but the package 'color.sty' is not loaded}%
    \renewcommand\color[2][]{}%
  }%
  \providecommand\transparent[1]{%
    \errmessage{(Inkscape) Transparency is used (non-zero) for the text in Inkscape, but the package 'transparent.sty' is not loaded}%
    \renewcommand\transparent[1]{}%
  }%
  \providecommand\rotatebox[2]{#2}%
  \newcommand*\fsize{\dimexpr\f@size pt\relax}%
  \newcommand*\lineheight[1]{\fontsize{\fsize}{#1\fsize}\selectfont}%
  \ifx\svgwidth\undefined%
    \setlength{\unitlength}{699.11486047bp}%
    \ifx\svgscale\undefined%
      \relax%
    \else%
      \setlength{\unitlength}{\unitlength * \real{\svgscale}}%
    \fi%
  \else%
    \setlength{\unitlength}{\svgwidth}%
  \fi%
  \global\let\svgwidth\undefined%
  \global\let\svgscale\undefined%
  \makeatother%
  \begin{picture}(1,0.65037959)%
    \lineheight{1}%
    \setlength\tabcolsep{0pt}%
    \put(0,0){\includegraphics[width=\unitlength,page=1]{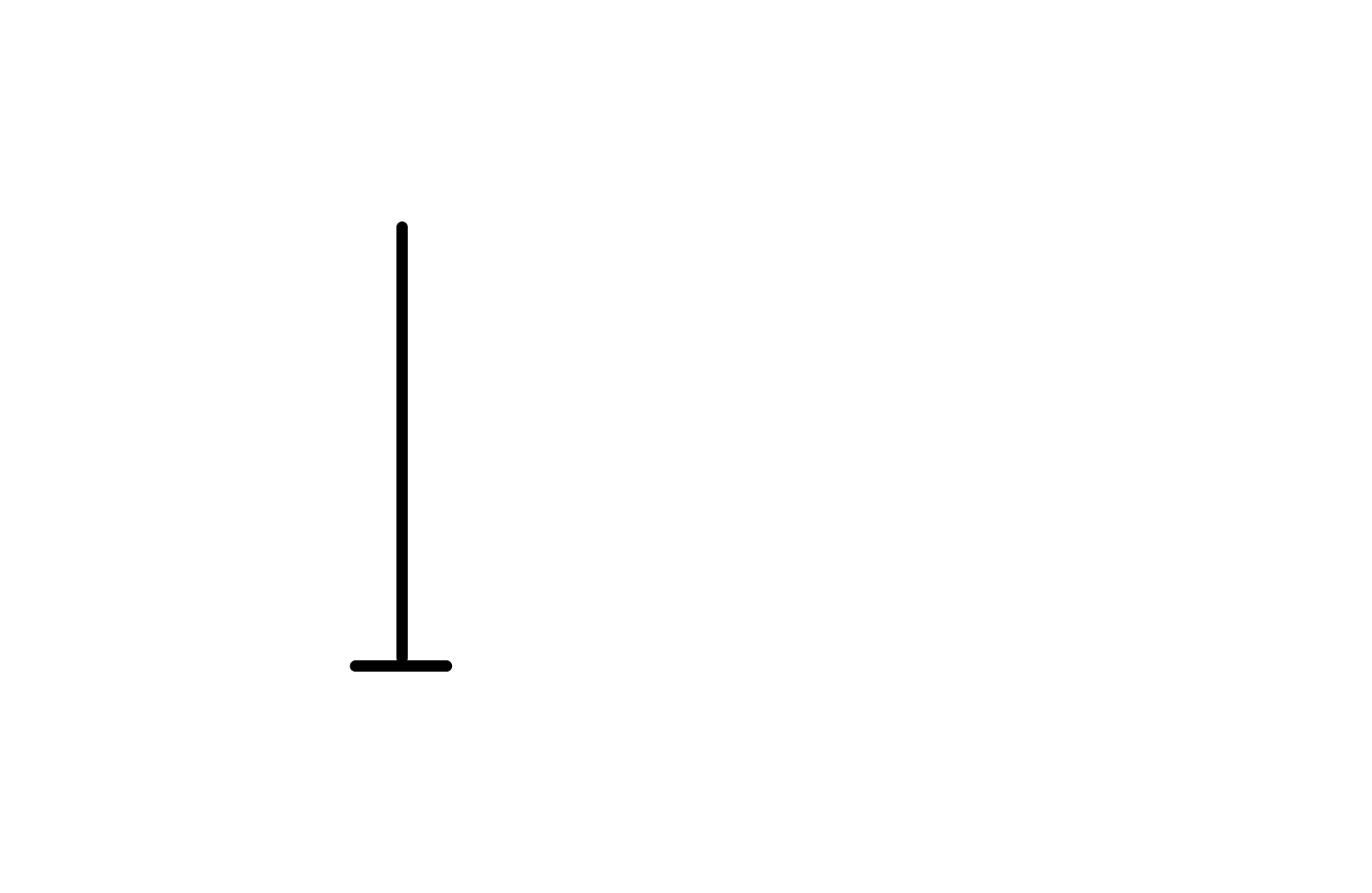}}%
    \put(0.46524565,0.3124238){\color[rgb]{0,0,0}\makebox(0,0)[lt]{\lineheight{1.25}\smash{\begin{tabular}[t]{l}$=$\end{tabular}}}}%
    \put(0,0){\includegraphics[width=\unitlength,page=2]{bialgebra4.pdf}}%
  \end{picture}%
\endgroup%

  \caption{}\label{fig bialgebra4}
\end{subfigure} 
   \hfill
\begin{subfigure}{0.12\linewidth}
\phantom{.}
\end{subfigure} 

\caption{ The bialgebra axioms are represented graphically. (\subref{fig associativity}), (\subref{fig unitality}) represent \eqref{eq algebra axioms}; (\subref{fig coassociativity}), (\subref{fig counitality}) represent \eqref{eq coalgebra axioms}, and  (\subref{fig bialgebra1}),(\subref{fig bialgebra2}),(\subref{fig bialgebra3}),(\subref{fig bialgebra4}) represent \eqref{eq bialgebra axioms}. Note that the right-hand side of (\subref{fig bialgebra4}) stands for the empty diagram.}
\label{fig bialgebra axioms}
\end{figure}

It will be convenient to depict the previous maps as some planar diagrams,
\begin{center}
\begin{tikzpicture}
\draw (0, 0) node[inner sep=0] {\includegraphics[width=1cm]{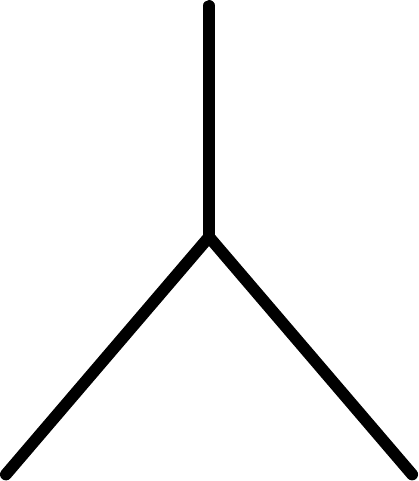}};
\draw (-1, 0) node {$\mu=$};
\draw (1.5, 0) node {, \quad $\eta=$};
\draw (2.5, 0) node[inner sep=0] {\includegraphics[width=1cm]{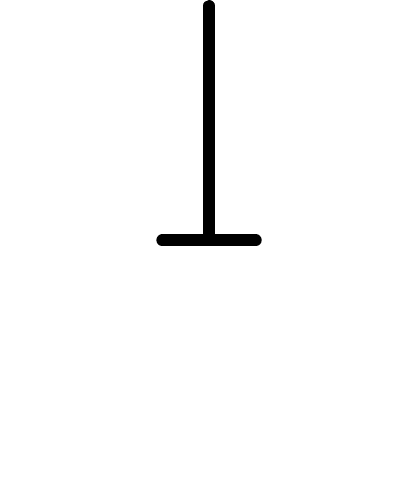}};
\draw (4, 0) node {, \quad $\Delta=$};
\draw (5.5, 0) node[inner sep=0] {\includegraphics[width=1cm]{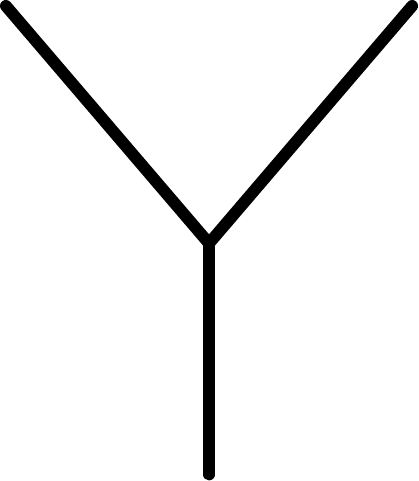}};
\draw (7, 0) node {, \quad $\varepsilon=$};
\draw (8, 0) node[inner sep=0] {\includegraphics[width=1cm]{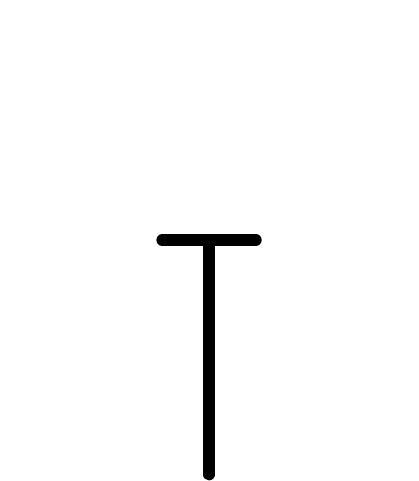}};
\end{tikzpicture}
\end{center}
whereas the symmetry will be represented as
\begin{center}
\begin{tikzpicture}
\draw (-1.5, 0) node {\quad $P_{A,A}=$};
\draw (0, 0) node[inner sep=0] {\includegraphics[width=1cm]{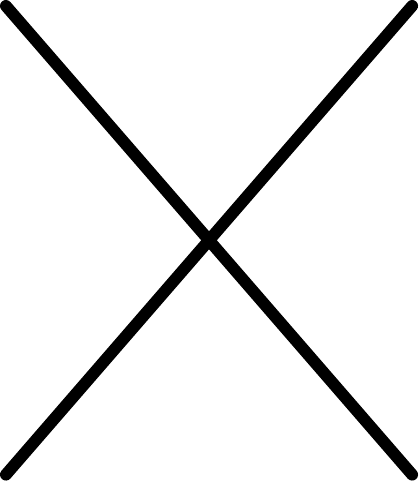}};
\end{tikzpicture}
\end{center}
By doing this, we can represent the bialgebra axioms graphically as in Figure \ref{fig bialgebra axioms}.

In the following we will adhere to Habiro's notation \cite{habiro2016}: for $k \geq 0$ we define the iterated multiplication and comultiplication $$\mu^{[k]}: A^{\otimes k} \to A \qquad , \qquad \Delta^{[k]} : A \to A^{\otimes k}$$ as follows:
\begin{gather*}
\mu^{[0]} :=  \eta , \quad \mu^{[1]} := \id \quad , \quad \mu^{[k]} := \mu \circ (\mu^{[k-1]} \otimes \id), \quad k \geq 2\\
\Delta^{[0]} :=  \varepsilon , \quad \Delta^{[1]} := \id \quad , \quad \Delta^{[k]} :=   (\Delta^{[k-1]} \otimes \id)\circ \Delta, \quad k \geq 2
\end{gather*}

For integers $p >0$ and $k_1, \ldots , k_p \geq 0$ let us write $$\mu^{[k_1, \ldots , k_p]} := \mu^{[k_1]} \otimes \cdots \otimes \mu^{[ k_p]} \qquad , \qquad \Delta^{[k_1, \ldots , k_p]} := \Delta^{[k_1]} \otimes \cdots \otimes \Delta^{[ k_p]}. $$
The associativity and coassociativity axioms imply the following generalised relations:
$$\mu^{[p]} \circ \mu^{[k_1, \ldots , k_p]} = \mu^{[k_1 + \cdots + k_p]} \qquad , \qquad   \Delta^{[k_1, \ldots , k_p]} \circ \Delta^{[p]} = \Delta^{[k_1+  \cdots + k_p]}. $$

\subsection{PROPs and their algebras}\label{subsec PROPs}

A \textit{(classical, monocoloured) PROP}  is a symmetric monoidal category  monoidally generated by a single object.  By MacLane's coherence theorem, every PROP is monoidally equivalent to another that has the set of non-negative integers as objects and the monoidal product is given by the sum.

PROPs are useful to study structural morphisms of large classes of algebraic structures. Let us review here a general construction to construct PROPs from generators and relations \cite{zanasi}. 

A \textit{symmetric monoidal theory} is a pair of sets $(\mathcal{G}, \mathcal{E})$,  where the \textit{signature} $\mathcal{G}$ contains \textit{generators}, a set of formal arrows $g: n \to m$ between non-negative integers, which includes two special arrows $\id: 1 \to 1$ and $P: 2 \to 2$. The set of \textit{$\mathcal{G}$-terms} is obtained by formally combining the generators with two operations ``composition'' $\circ$ and ``monoidal product'' $\otimes$ as follows: given $\mathcal{G}$-terms $f: n \to m$, $g: m \to s$ and $h: p \to q$, we construct new $\mathcal{G}$-terms $g \circ f: n \to s$ and $f  \otimes h: n+p  \to m+q$. The set $\mathcal{E}$ contains \textit{equations}, pairs of $\mathcal{G}$-terms $(t,t':n \to m )$ with the same arity and coarity.

Any symmetric monoidal theory $(\mathcal{G}, \mathcal{E})$ freely generates a PROP $\mathsf{P}=\mathsf{P}_{(\mathcal{G}, \mathcal{E})}$ by letting $\hom{\mathsf{P}}{n}{m}$ be the set of $\mathcal{G}$-terms $n \to m$ modulo the least congruence relation (with respect to the composition and the monoidal product)  that contains the laws of strict symmetric categories and the equations $t=t'$ for every pair $(t,t') \in \mathcal{E}$. In other words, $\mathsf{P}$ is the strict symmetric monoidal category monoidally generated by the arrows of $\mathcal{G}$ subject to the relations given by $\mathcal{E}$.

Whereas a PROP encodes algebraic operations abstractly, an algebra over a PROP is an evaluation of these operations on a concrete object. More precisely, a \textit{$\mathsf{P}$-algebra} in a symmetric monoidal category $\C$ is a strong monoidal functor ${A: \mathsf{P} \to \C}$. The class of $\mathsf{P}$-algebras on $\C$  forms a category $\mathsf{Alg}_{\mathsf{P}} (\C)$ in the obvious way, where morphisms are monoidal natural transformations between functors $ \mathsf{P} \to \C$.

The following observation will be useful:

\begin{lemma}\label{lemma yoneda}
Let $\mathsf{P}$, $\mathsf{P}'$ be PROPs. If there is a natural bijection $$\mathsf{Alg}_{\mathsf{P}} (\C) \toiso \mathsf{Alg}_{\mathsf{P}'} (\C)$$ for any symmetric monoidal category $\C$, then $\mathsf{P}$ and $\mathsf{P}'$ are monoidally equivalent.
\end{lemma}
\begin{proof}
This is a direct consequence of the Yoneda lemma.
\end{proof}

Let us illustrate the above construction with a few examples:

\begin{example}
Let $\mathcal{G}$ consist of two arrows $\mu: 2\to 1$ and $\eta: 0 \to 1$, that we can think of as the planar diagrams depicted in the previous section. If $\mathcal{E}$ is the set given by the two equations depicted in Figure \ref{fig associativity} and \ref{fig unitality}, then the resulting PROP is denoted by $\mathsf{A}$ and $\mathsf{Alg}_{\mathsf{A}} (\C)$ is equivalent to the category of algebras in $\C$. In this case, one says that $\mathsf{A}$ \textit{is a PROP for algebras}. If the relation $P \circ  \mu =\mu$ is added, then the corresponding category $\mathsf{ComA}$ is a PROP for commutative algebras.

A similar note can be made for the categories of coalgebra or bialgebras, perhaps with the (co)commutative axiom included. If $\mathcal{G}$ contains $\mu:2\to 1$, $\eta: 0\to 1$, $\Delta: 1\to 2$ and $\varepsilon: 1 \to 0$ and $\mathcal{E}$ consists of the equations given in Figure \ref{fig bialgebra axioms}, then the resulting PROP for bialgebras will be denoted as $\mathsf{B}$. If the equation $P \circ  \mu =\mu$ is added, then we will denote the resulting PROP for commutative bialgebras as $\mathsf{ComB}$.
\end{example}

An important drawback of the previous construction is that it  might seem ad-hoc or artificial, being produced by generators and relations. A sensible question is the following: 

\begin{question}
Is it possible to describe the PROPs for ((co)commutative) algebras, coalgebras or bialgebras without using generators and relations?
\end{question}

The answer for ((co)commutative) (co)algebras is quite simple. Before addressing this case let us introduce some notation: for $n \geq 0$, let $\Sigma_n$ be the symmetric group of order $n$ (if $n=0$ we set $\Sigma_0 = \emptyset$).  Given an object $X \in \C$, any element $\sigma \in \Sigma_n$ induces an arrow $$P_\sigma: X^{\otimes n} \to X^{\otimes n}$$ determined by $P_{(i,i+1)}=\id^{\otimes (i-1)} \otimes P_{X,X} \otimes \id^{n-i-1}$ and the property that the passage $\Sigma_n \to \hom{\C}{X^{\otimes n}}{X^{\otimes n}}$ is a monoid homomorphism.

\begin{example}[essentially \cite{pirashvili}]\label{ex fSet A}
Let $\mathsf{fSet}$ be the  category of finite sets. Up to monoidal equivalence, we can view it as a PROP with the ordinals $[n] = \{ 1< \cdots < n  \}$, $n \geq0$, as objects and the sum of ordinals (disjoint union)  as the monoidal product ($[0]$ is by definition the empty set).

For any $k \geq 0$, there is a unique map $[k] \to [1]$, in the same way that there is a unique map $\mu^{[k]}: A^{\otimes k} \to A$ for a commutative algebra $A$ in $\C$.\footnote{Homotopy theorists will recognise the commutative operad $\mathsf{Com}(n)= *$ here.} This suggests the following: given a strong monoidal functor $F: \mathsf{fSet} \to \C$, let $X:= F([1])$ and let $\mu^{[k]} := F([k] \to [1])$. This determines a structure of a commutative algebra object on $X$. Conversely, if $A\in \C$ is a commutative algebra, we define a strong monoidal functor $F: \mathsf{fSet} \to \C$ as follows: first set $F([1]):=A$. For a map of sets $f:[n] \to [m]$, let $k_i := \#f^{-1}(i) \geq 0$ for every $i \in [m]$. If $\sigma \in \Sigma_n$ is any permutation that maps the subset $\{k_1 + \cdots + k_{i-1} +1 , \ldots , k_1 + \cdots + k_i \}$ to $f^{-1}(i)$, then define $$Ff := \mu^{[k_1, \ldots , k_m]} \circ  P_\sigma,$$ where $P_\sigma: A^{\otimes n} \to A^{\otimes n}$ is as above. The commutativity of $A$ ensures that any choice of $\sigma \in \Sigma_n$ yields the same arrow. This shows, by \ref{lemma yoneda}, that $\mathsf{fSet}$ is equivalent to the PROP $\mathsf{ComA}$ for commutative algebras.

Let us illustrate the previous construction in the case $\C = \mathsf{Vect_k}$, the category of vector spaces over a field $k$. If $A$ is an algebra and $f:[n] \to [m]$ then $$(Ff)(a_1 \otimes \cdots \otimes a_n)= \left( \prod_{i \in f^{-1}(1)} a_i \right) \otimes \cdots \otimes \left( \prod_{i \in f^{-1}(m)} a_i \right),$$ where an empty product stands for $1 \in A$.

By reversing arrows, we have that $\mathsf{fSet}^{op}$ is equivalent to the PROP $\mathsf{CocomC}$ for co\-co\-mmu\-ta\-tive coalgebras.
\end{example}

\begin{example}[\cite{pirashvili}]\label{ex pirashvili ordered}
Now suppose we want to refine the previous example to obtain the PROP $\mathsf{A}$ for (non-commutative) algebras. If  the order of multiplication matters, we are forced to impose orders on each fibre $f^{-1}(i)$, which in turn determines a unique permutation. 

More concretely, define  $\widetilde{\mathsf{fSet}}$ as the following category: its objects are the ordinals $[n]$, $n \geq 0$. An arrow $f: [n] \to [m]$ is the data of a map of sets $f: [n] \to [m]$ together with a total order of $f^{-1}(i)$ for every $i \in [m]$. Given arrows $f:[n] \to [m]$ and $g: [m] \to [p]$, the composite $g \circ f$ is defined as the composite of the underlying maps of sets together with the following order on fibres: for $i \in [p]$, let 
\begin{equation}\label{eq composite in fSettilde}
(g \circ f)^{-1}(i) := f^{-1}(j_1) \amalg \cdots \amalg f^{-1}(j_r)
\end{equation}
where $g^{-1}(i)= \{j_1 < \cdots < j_r \}$ and the disjoint union carries the order imposed by the writing from left to right. The disjoint union makes it a symmetric monoidal category in the obvious way. Reasoning as in the previous example, we see that given a map $f$ there is a unique $\sigma \in \Sigma_n$ mapping the ordered set $\{k_1 + \cdots + k_{i-1} +1 < \cdots < k_1 + \cdots + k_i \}$ to $f^{-1}(i)$ in an order-preserving way\footnote{This is another way of saying that the associative operad $\mathsf{Ass}(n)$ equals $\Sigma_n$.}, and therefore $\widetilde{\mathsf{fSet}}$ is equivalent to the PROP for algebras.

As a concrete instance of this construction, let
\begin{align*}
f:[5]\to [4] \quad &, \quad f^{-1}(1)= \{ 4<2 \} \quad , \quad f^{-1}(3)= \{ 1<3 \} \quad , \quad f^{-1}(4)= \{ 5 \},\\
g:[4]\to [2] \quad &, \quad g^{-1}(1)= \{ 4<3 \} \quad , \quad g^{-1}(2)= \{ 1<2 \} .
\end{align*}
Then the composite is determined by $$(g \circ f)^{-1}(1) = \{ 5<1<3 \} \quad , \quad (g \circ f)^{-1}(2) = \{ 4<2\}.$$
The order on fibres that $f$ determines yields a unique permutation $\sigma = (143) \in \Sigma_5$ such that $Ff = \mu^{[2, 0,2,1]} \circ  P_\sigma$. Likewise, $g$ determines a unique $\tau = (1423) \in \Sigma_4$ such that $Fg = \mu^{[2,2]} \circ P_\tau$. If $\C= \mathsf{Vect_k}$, then $(Ff)(a_1 \otimes \cdots \otimes a_5) = a_4a_2 \otimes 1 \otimes a_1a_3\otimes a_5$ and $(Fg)(b_1 \otimes \cdots \otimes b_4) = b_4b_3 \otimes b_1b_2 $.
\end{example}



What we learned from this example is that by adding extra data to the PROP for commutative algebras, we were able to ``fix the order of multiplication'' and obtain a PROP for algebras. We also see that the data of a total order on fibres is equivalent to the datum of a permutation, since for $f: [n] \to [m]$ and $\sigma$ as before $$f^{-1}(i)= \{\sigma (k_1 + \cdots + k_{i-1} +1) < \cdots < \sigma (k_1 + \cdots + k_i) \}.$$ We can then replace the data of one by the data of the other.

\begin{construction}
Let $\alpha \in \Sigma_m$ and let $k_1, \ldots, k_m \geq 0$ be a sequence of $m$ non-negative integers. If $n= \sum_i k_i$, then $\alpha$ induces a permutation  $\langle \alpha \rangle_{k_1, \ldots, k_m} \in \Sigma_n$ which permutes the $m$ consecutive blocks of $[n]$ with sizes $ k_1, \ldots, k_m$. More precisely, for $i \in [m]$, define $\langle \alpha \rangle_{k_1, \ldots, k_m}$ as  the unique permutation that sends $\{k_{\alpha(1)} + \cdots + k_{\alpha(i-1)}+1, \ldots , k_{\alpha(1)} + \cdots + k_{\alpha(i)} \}$ to $\{k_1 + \cdots + k_{\alpha(i)-1}+1, \ldots , k_1 + \cdots + k_{\alpha(i)} \}$ in a order-preserving way.


A schematic of this construction  is  shown in Figure \ref{fig langle rangle constr}. It is immediate to check that $$\langle \alpha \rangle_{1, \overset{m}{\ldots}, 1} = \alpha \qquad , \qquad \langle \id_m \rangle_{k_1, \ldots, k_m} = \id_n$$ for any $\alpha \in \Sigma_m$ and $k_1, \ldots , k_m \geq 0$.

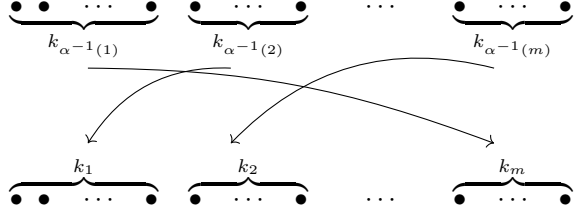
\begin{figure}[h]
\begin{tikzpicture}
\draw (0, 0) node {$\underbrace{ \bullet  \  \bullet \quad \cdots \quad \bullet}_{k_{\alpha^{-1}(1)}} \quad \underbrace{ \bullet  \quad \cdots \quad  \bullet }_{k_{\alpha^{-1} (2)}}  \qquad \cdots \qquad \underbrace{ \bullet  \quad \cdots \quad  \bullet }_{k_{\alpha^{-1} (m)}} $};
\draw (0, -2) node {$\overbrace{ \bullet  \   \bullet  \quad \cdots \quad  \bullet }^{k_1} \quad \overbrace{ \bullet  \quad \cdots \quad  \bullet }^{k_2} \qquad \cdots \qquad \overbrace{ \bullet  \quad \cdots \quad  \bullet }^{k_m} $};

\draw [->] (-0.8,-0.5) to [bend right] (-2.7,-1.5);
\draw [->] (2.7,-0.5) to [bend right] (-0.8,-1.5);
\draw [->] (-2.7,-0.5) to [bend left=10] (2.7,-1.5);

\end{tikzpicture}
\caption{Illustration of the construction of $\langle \alpha \rangle_{k_1, \ldots, k_m}$. If we label the bullets of the upper row with the sequence $1, \ldots, n$, and we take each of the blocks to the row below according to $\alpha$, that is, identifying in the two rows blocks with the same index $k_{i}$, then the resulting labelling of the lower row contains, from left to right, the values $\langle \alpha \rangle_{k_1, \ldots, k_m} (1), \ldots , \langle \alpha \rangle_{k_1, \ldots, k_m}(n)$.}
\label{fig langle rangle constr}
\end{figure}
\end{construction}

We also introduce the following notation: given two permutations $\alpha \in \Sigma_n$, ${\beta \in \Sigma_m}$, we denote by $\alpha \otimes \beta \in \Sigma_{n+m}$ the \textit{block product} of $\alpha$ and $\beta$, that is, the image of the the pair $(\alpha , \beta)$ under the ``block inclusion'' $$\Sigma_n \times \Sigma_m \hooklongrightarrow \Sigma_{n+m}.$$ Concretely, $$(\alpha \otimes \beta) (i) = \begin{cases}
\alpha (i), & i=1, \ldots ,n\\
\beta (i-n)+n, & i=n+1, \ldots , n+m.
\end{cases} $$

\begin{definition}
Let $\widehat{\mathsf{fSet}}$ be the following category: its objects are the ordinals. An arrow is a pair $(f, \sigma)$ where $f$ is a map of sets $f: [n] \to [m]$  and $\sigma \in  \Sigma_n$ is a permutation.  Given pairs $(f, \sigma)$ and $(g, \tau)$ with $f$ as before and $g: [m] \to [p]$, define the composite as $$(g, \tau) \circ (f, \sigma):= (g \circ f,   \sigma \circ \langle  \tau \rangle_{k_1, \ldots , k_m}) $$ where $k_i := \#f^{-1} (i)$ for every $i \in [m]$. The identity is the pair $(\id, \id)$. Furthermore there is a monoidal product given by $$(f_1, \sigma_1) \otimes (f_2, \sigma_2) := (f_1 \amalg f_2, \sigma_1 \otimes \sigma_2).$$ 
\end{definition}

\begin{proposition}\label{prop fSethat fSettilde}
The categories $\widehat{\mathsf{fSet}}$ and $\widetilde{\mathsf{fSet}}$ are monoidally equivalent, so $\widehat{\mathsf{fSet}}$ is also equivalent to $\mathsf{A}$.
\end{proposition}
\begin{proof}
We have already described how giving an order on fibres amounts to giving a permutation. Hence to see that this defines a functor $\widetilde{\mathsf{fSet}} \to \widehat{\mathsf{fSet}} $ all is left to check is that the correspondence preserves the composition.

Let $f:[n] \to [m]$, $g:[m]\to [p]$ be composable arrows in $\widetilde{\mathsf{fSet}}$ with $\sigma \in \Sigma_n$, $\tau \in \Sigma_m$ the permutations that arise from the order of the fibres as explained in \ref{ex pirashvili ordered}, respectively. If $\delta$ is the permutation associated to the composite, then it is   determined by the following equality of ordered sets:
\begin{align*}
\{ \delta(1) < \cdots < \delta(n)  \} &= (g \circ f)^{-1}(1) \amalg \cdots \amalg (g \circ f)^{-1}(p)\\
&= \coprod_{j \in g^{-1}(1) \amalg \cdots \amalg g^{-1}(p)}    f^{-1}(j)\\
&=  \coprod_{j \in \{ \tau(1) < \cdots < \tau (m)  \}} \{\sigma (k_1 + \cdots + k_{j-1} +1) < \cdots < \sigma (k_1 + \cdots + k_j) \}.
\end{align*}
This expresses that $\delta$ is given by a block permutation of the set $\{ \sigma(1) < \cdots < \sigma(n)  \}$, where the blocks have sizes $k_1 , \ldots , k_m$ and the order of the block permutation is determined by $\tau$. This exactly says that $\delta =\sigma \circ \langle  \tau \rangle_{k_1, \ldots , k_m}$.
\end{proof}


\subsection{Finitely generated free monoids}

Recall that there is a free - forgetful adjunction
$$\begin{tikzcd}[column sep={4em,between origins}]
\mathsf{Set}
  \arrow[rr, bend left, swap, "F"' pos=0.5]
& \bot &   \mathsf{Mon}
  \arrow[ll, bend left, swap, "U"' pos=0.5]
\end{tikzcd}$$
between the category of sets and the category of monoids, where $F(X):= \coprod_{n \geq 0} X^{\times n}$ is given by ``words on the alphabet $X$''. We denote by $\mathsf{fgFMon}$ the category of finitely-generated free monoids, that is, the full subcategory of $\mathsf{Mon}$ on the objects $F([n])$, $n \geq 0$.

For the sake of clarity let us label $[n]=\{x_1 < \cdots < x_n  \}$ with ``$x$'s'' or any other letter so that $$F(x_1, \ldots , x_n) := F([n])= F(\{ x_1, \ldots , x_n  \}).$$ Observe that by the adjunction,
\begin{align*}
\hom{\mathsf{Mon}}{F(x_1, \ldots , x_n)}{F(y_1, \ldots , y_m)} &\cong \hom{\mathsf{Set}}{\{ x_1, \ldots , x_n  \}}{UF(y_1, \ldots , y_m)} \\
&\cong \prod_{i=1}^n \hom{\mathsf{Set}}{\{ x_i  \}}{UF(y_1, \ldots , y_m)}\\
&\cong \prod_{i=1}^n UF(y_1, \ldots , y_m),
\end{align*}
so any monoid map between free monoids is completely determined by a tuple of words (the images of the elements $x_i$).

The category  $\mathsf{fgFMon}$ has finite coproducts given by the free product of monoids, $$F(x_1, \ldots , x_n) * F(y_1, \ldots , y_m) \cong F(x_1, \ldots , x_n,y_1, \ldots , y_m ),$$ which turns it into  a cocartesian monoidal category, where the initial object $\mathbf{1}:= F([0])$ (the trivial monoid) serves as the unit of the monoidal product.

\begin{example}[\cite{habiro2016}]\label{ex ComB = fgFMon}
The object $F(x)= F([1]) \in \mathsf{fgFMon}$ is a bialgebra object: the multiplication and multiplication are given by
\begin{align*}
&\mu: F(x,y) \to F(z) \qquad , \qquad x,y \mapsto z,\\
&\Delta: F(x) \to F(y,z) \qquad , \qquad x \mapsto yz,
\end{align*}
and the unit and counit are the unique maps $\eta: \mathbf{1} \to F(x)$ and $\varepsilon:  F(x)\to \mathbf{1}$. It is straightforward to check that the bialgebra axioms \eqref{eq algebra axioms} -- \eqref{eq bialgebra axioms} hold. For instance, both sides of the associativity axiom from Figure \ref{fig bialgebra axioms}(\subref{fig  associativity}) are the map $x,y,z \mapsto x$; from Figure \ref{fig bialgebra axioms}(\subref{fig  coassociativity})  they are the map $x \mapsto xyz$; and from Figure \ref{fig bialgebra axioms}(\subref{fig bialgebra1}) they are the map $x, y \mapsto ab$. The rest of axioms can be similarly checked.

Therefore there is a functor $$T: \mathsf{ComB} \to \mathsf{fgFMon}.$$
It is a well-know result that the above functor is a monoidal equivalence. Habiro \cite{habiro2016} showed this directly\footnote{Actually, Habiro showed the more general result that the PROP for commutative Hopf algebras is equivalent to the category of finitely generated free groups, but his construction easily restricts to the bialgebra case.}, but it can also be proven using Lawvere theories \cite{pirashvili}.
\end{example}

To finish this section we would like to explain how the equivalences $\mathsf{ComA} \simeq \mathsf{fSet}$ and  $\mathsf{ComB} \simeq \mathsf{fgFMon}$ are related. Firstly, observe that there is a canonical faithful functor
\begin{equation}\label{eq ComA Comb}
\mathsf{ComA} \hooklongrightarrow \mathsf{ComB}
\end{equation}
which arises from the inclusion of generators and relations between the symmetric monoidal theories freely generating both PROPs.

\begin{lemma}\label{lem fSet fgFMon}
The free monoid functor $F: \mathsf{fSet} \to \mathsf{fgFMon}$ makes the following diagram commutative,
$$\begin{tikzcd}
\mathsf{ComA} \dar{\simeq} \rar[hook] &  \mathsf{ComB}  \dar{\simeq}\\
\mathsf{fSet} \rar{F} & \mathsf{fgFMon}
\end{tikzcd}   $$
where the upper arrow is the natural inclusion.
\end{lemma}
\begin{proof}
This is immediate: the generator $\mu^{[k]}: k \to 1$ of $\mathsf{ComA}$, $k \geq 0$, is mapped under the two possible composites to the monoid map $F(x_1, \ldots , x_k) \to F(y)$ determined by $x_i  \mapsto y$ for all $1 \leq i \leq k$, the iterated multiplication in  $\mathsf{fgFMon}$.
\end{proof}

\begin{remark}
The non-existence of a right-adjoint for the free functor $F: \mathsf{fSet} \to \mathsf{Mon}$ can be seen as the failure of the inclusion $\mathsf{ComA}  \hooklongrightarrow \mathsf{ComB}$ to have a right-adjoint. 
\end{remark}


\section{A PROP for bialgebras}\label{sec A PROP for bialgebras}

The main question we want to address is whether it is possible to obtain a PROP for (non)commutative bialgebras out of   $\mathsf{fgFMon}$ by adding extra data, as we did in the algebra case. We would like to start with some motivating examples:

\begin{example}\label{ex motivating 1}
Let $A$ be a commutative algebra in $\mathsf{Vect}_k$, and let $\mathsf{fgFMon}\to \mathsf{Vect}_k$ be the functor making the diagram
$$\begin{tikzcd}
\mathsf{ComB} \arrow{rr}{\simeq} \drar[swap]{A} & & \mathsf{fgFMon} \dlar\\
&\mathsf{Vect}_k &
\end{tikzcd}$$
commutative (up to natural isomorphism), so that we can view $A$ as a functor $A: \mathsf{fgFMon}\to \mathsf{Vect}_k$. Let us investigate this point of view: the non-cocommutativity of the comultiplication $\Delta: A \to A \otimes A$, $$\sum_{(a)} a_{(1)} \otimes a_{(2)} \neq \sum_{(a)} a_{(2)} \otimes a_{(1)},$$  can be seen as a consequence of the non-equality of the monoid maps $F(x) \to F(y,z)$ given by $(x \mapsto yz)$ and $(x \mapsto zy)$. 

Similarly, the commutativity of the multiplication $\mu: A \otimes A \to A$, $ab = ba$, can be seen as a consequence of the commutativity of the following diagram:
$$\begin{tikzcd}
F(x,y) \arrow{rr}{P} \drar[swap]{\mu} & & F(x,y) \dlar{\mu}\\
&F(z) &
\end{tikzcd}$$
where $P=(x\mapsto y, y \mapsto x)$. So producing non-commutative bialgebras amounts to being able to distinguish $\mu$ and $\mu \circ P $. One way to achieve this is as follows: note that the ordinals $[n] = \{1 < \cdots <n \}$ carry a canonical order, which induces an order\footnote{Usually we write $\{x_1, x_2, x_3, \ldots \} = \{a, b,c, \ldots \}$. In the latter case the order will always be alphabetical.} on the subset  $\{ x_1, \ldots , x_n \} \subset F(x_1, \ldots , x_n)$. In particular, there is a ``canonically ordered'' word $x_1 \cdots x_n \in F(x_1, \ldots , x_n)$. For the multiplication ${\mu: F(x,y) \to F(z)}$, this element $xy$ maps to $z^2$. This corresponds to $a\otimes b \mapsto ab=ba$ one we pass to the algebra. If we are interested in fixing the order of multiplication, we could associate the permutation $\id_2 \in \Sigma_2$ to the first option $a \otimes b \mapsto ab$ and $(12) \in \Sigma_2$ to the second one, $a \otimes b \mapsto ba$. At the level of the free monoids, this corresponds to associating $\id_2$ to $\mu$ and $(12)$ to $\mu \circ P$.
\end{example}

\begin{example}
Let $f: F(x) \to F(y)$ be a monoid map, which necessarily must be  $x \mapsto y^n$ for some $n \geq 0$. Given a strong monoidal functor  $A: \mathsf{fgFMon}\to \mathsf{Vect}_k$, this map will be sent to the map $A \to A$, $a \mapsto a_{(1)} \cdots a_{(n)}$ (one way to see this is that $f$ factors as $f= \mu^{[k]} \Delta^{[k]}$). Since $A$ is a commutative algebra, then $ a_{(1)} \cdots a_{(n)} = a_{(\sigma (1))} \cdots a_{(\sigma (n))}$ for any $\sigma \in \Sigma_n$. If we want to lift this correspondence to the non-commutative case, we are forced to fix the indeterminacy that the assignment $x \mapsto y^n$ carries for the order of multiplication. If we fix  $\sigma \in \Sigma_n$, then we also fix the order of multiplication.
\end{example}

The above two examples suggest that to extend $\mathsf{fgFMon}$ to a PROP for (non-commutative) bialgebras, we have to add the extra data of some permutations that will determine the order of multiplication once we evaluate in a algebra.

\begin{notation}\label{notation}
Let $A$ be a bialgebra in $\mathsf{Vect}_k$ and let $f: n \to m$ be an arrow in $\mathsf{B}$. If $A: \mathsf{B} \to \mathsf{Vect}_k$ as a functor, then it is possible to read off $Ff$ (and hence $f$) from the expression of $(Af)(a_1 \otimes \cdots \otimes a_n) \in A^{\otimes m}$. For instance, if $(a \otimes b \mapsto \varepsilon (a_{(2)}) b_{(1)} \otimes 1 \otimes b_{(2)}a_{(1)}  )$ then $Af$ must be
$$Af = (\varepsilon \otimes \id \otimes \eta \otimes \mu) \circ P_{(1234)} \circ (\Delta \otimes \Delta)  $$ (obviously $f$ has the same expression). Therefore,  for a algebra object $A$ of a symmetric monoidal category $\C$, we will adopt the convention of writing maps $Af$ with the symbols that are used for $k$-algebras. When convenient, we will denote a map in $\C$ as eg $(a \otimes b \mapsto \varepsilon (a_{(2)}) b_{(1)} \otimes 1 \otimes b_{(2)}a_{(1)}  )$ keeping in mind that there is no underlying  set in $A$ and this is purely a notational convenience.
\end{notation}

\subsection{A unique normal form in $\mathsf{B}$}

On our way to obtain an equivalent PROP for bialgebras it will become essential to have control over the arrows we can find in $\mathsf{B}$. The main result of this section ensures that any such an arrow has a unique ``normal form''. The main tool to show this will be the so-called ``diamond lemma''. Let us recall the setup.

Let $S$ be a set. A \textit{reduction} on $S$ is a strictly antisymmetric relation on $S$, that is, a subset $R \subset S \times S$ with the property that  $(x,y) \in R$ implies $(y,x) \not\in R$, and we will simply  denote this relation  as $\rightarrow$. A \textit{reduction chain} on $S$ is a sequence $x_1 \rightarrow x_2 \rightarrow  x_3 \rightarrow \cdots$. We say that the relation $\rightarrow$ satisfies the \textit{descending chain condition} if every reduction chain is finite. An element $x \in S$ is said to be in \textit{normal form} if it is minimal with respect to $\rightarrow$, that is, if there is no $y \in S$ such that $x \rightarrow y$. A subset $T \subset S$ is said to be \textit{connected} if whenever $x \in T$ and $y \in S$ with $x \rightarrow y$ or $y \rightarrow x$, then $y \in T$.  Given a reduction  $\rightarrow$ on $S$, its \textit{transitive-reflexive closure}  is the order relation $\twoheadrightarrow$ defined as $x \twoheadrightarrow y$ if and only if there is a reduction chain connecting $x$ and $y$. More formally, $x \twoheadrightarrow y$ if and only if there exist $x_0, \ldots , x_n \in S$ such that $x=x_0  \rightarrow  x_1 \rightarrow  \cdots \rightarrow  x_n=y$.  Moreover, if $x  \twoheadrightarrow y$ and $y$ is in normal form, then we say that $y$ \textit{is a normal form of} $x$. Finally, we say that $\rightarrow$ satisfies the \textit{diamond condition} if given $y \leftarrow x \rightarrow z$, there exists  $w \in S$ such that $y \twoheadrightarrow w$ and $z \twoheadrightarrow w$. 

\begin{lemma}[Diamond lemma]\label{lem diamond} 
Let $\rightarrow$ be a reduction on a non-empty set $S$. If $\rightarrow$ satisfies the descending chain condition and the diamond condition, then every element in $S$ has a unique normal form. 
\end{lemma}

We refer the reader to \cite[4.6]{barnatanveen} or \cite[4.47]{becker} for a proof. It readily follows from the lemma that all elements of a given connected subset have the same (unique) normal form.

\begin{theorem}\label{thm normal form}
Every morphism $f: n \to m$ in $\mathsf{B}$ factors in a unique way as $$f= \mu^{[q_1, \ldots , q_m]} \circ P_{\sigma} \circ \Delta^{[p_1, \ldots , p_n]}  $$ for some (unique) integers $p_1, \ldots , p_n, q_1, \ldots , q_m \geq 0$ such that $\sum_i p_i = \sum_i q_i =:s$, and a (unique) permutation $\sigma \in \Sigma_s$.
\end{theorem}

\begin{figure}[h]
 \centering
\def\svgwidth{0.45\textwidth}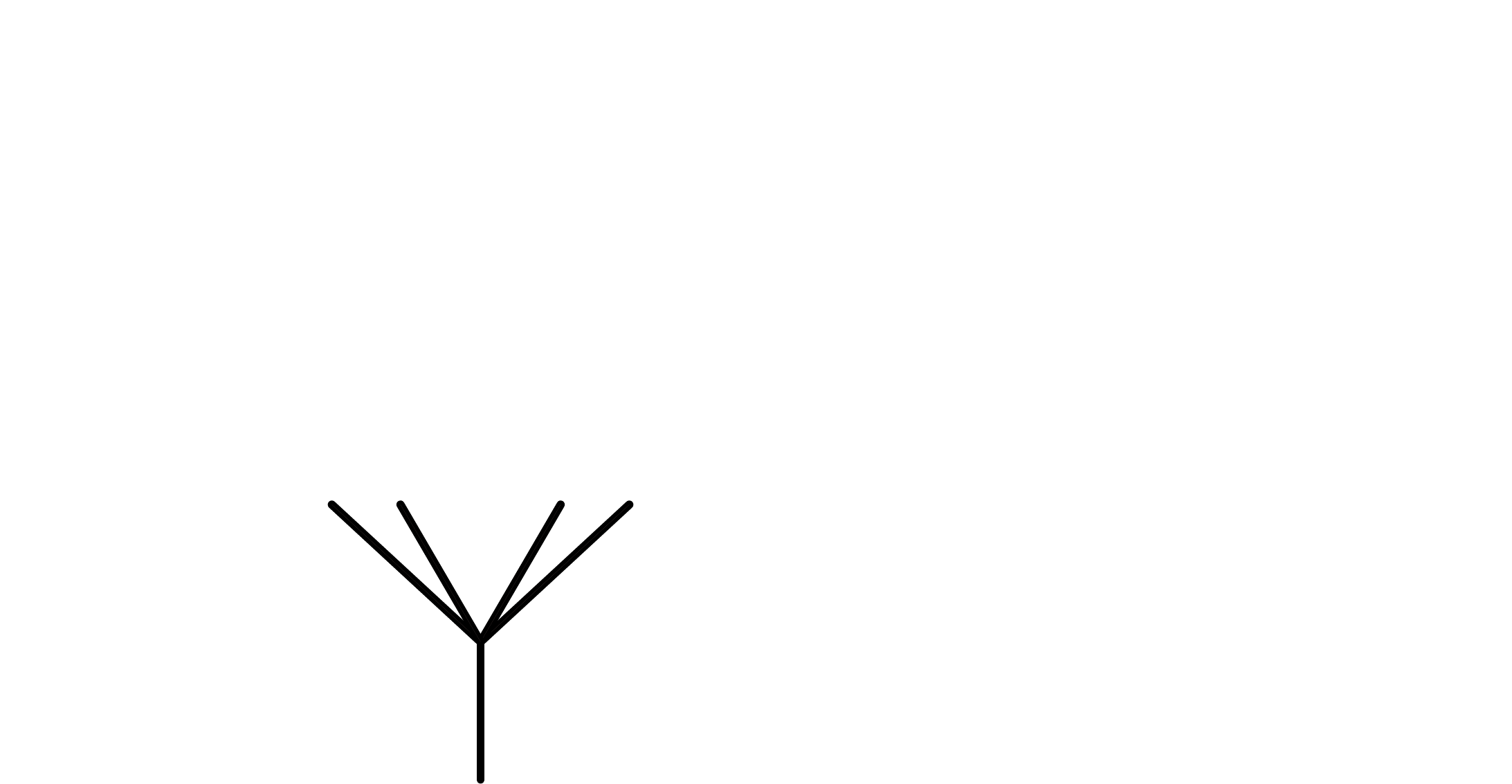
\caption{}
\label{fig factorisation}
\end{figure}

Graphically, we can represent this factorisation as depicted in Figure \ref{fig factorisation}.

The observation that ignites the proof is that in the relations that define the bialgebra PROP $\mathsf{B}$, depicted in Figure \ref{fig bialgebra axioms}, with the exception of (a) and (c) (since in the statement of the above theorem the iterated (co)multiplication appears), only one of the sides of the equations are in the normal form that the theorem suggests.

Let $\mathcal{B}$ be the PROP freely generated by the symmetric monoidal theory $(\mathcal{G},\mathcal{E})$ where the signature  $\mathcal{G}$ is given by arrows  $\mu:2\to 1$, $\eta: 0\to 1$, $\Delta: 1\to 2$ and $\varepsilon: 1 \to 0$, and  $\mathcal{E}$ contains the  equations appearing  in Figure \ref{fig bialgebra axioms}(a),(c) (the construction of  such PROPs was described in subsection \ref{subsec PROPs}). We are going to define a reduction on the set $$\mathrm{arr} (\mathcal{B}) := \coprod_{n,m \geq 0} \hom{\mathcal{B}}{n}{m}. $$

Let
\begin{align*}
\mu(\id \otimes \eta) \rightarrow \id \leftarrow \mu(\eta \otimes \id) \qquad &, \qquad  (\id \otimes \varepsilon)\Delta \rightarrow  \id \leftarrow  (\varepsilon \otimes \id)\Delta\\ 
\Delta \mu \rightarrow (\mu \otimes \mu)(\id \otimes P_{1,1} \otimes \id)(\Delta \otimes \Delta) \qquad &, \qquad  \Delta \eta \rightarrow \eta \otimes \eta\\
\varepsilon \mu \rightarrow \varepsilon \otimes \varepsilon \qquad &, \qquad \varepsilon \eta \rightarrow \id_{\mathds{1}}
\end{align*}
be the \textit{elementary reductions}, and extend $\rightarrow$ to $\mathrm{arr} (\mathcal{B})$ in way compatible with the composite and the monoidal product. More precisely, consider over $\mathrm{arr} (\mathcal{B})$ the smallest reduction $\rightarrow$ that contains the elementary reductions and that it is closed under composition and monoidal product in the sense that   if $f \rightarrow g$ is a reduction then so is $$h_4 \circ (h_2 \otimes f \otimes h_3) \circ h_1 \rightarrow h_4 \circ (h_2 \otimes g \otimes h_3) \circ h_1  $$
for any arrows $h_1, \ldots, h_4 \in \mathrm{arr} (\mathcal{B}) $ (so long as the composite makes sense).

\begin{lemma}\label{lem elements in normal form}
The subset of elements of $\mathrm{arr} (\mathcal{B})$ in normal form  with the reduction described above equals the set  $$ \{  \mu^{[q_1, \ldots , q_m]} \circ P_{\sigma} \circ \Delta^{[p_1, \ldots , p_n]} :  p_1, \ldots , p_n, q_1, \ldots , q_m \geq 0,  \sum_i p_i = \sum_i q_i =:s , \sigma \in \Sigma_s \}.$$
\end{lemma}
\begin{proof}
It is clear that the elements of this set are in normal form, since there is no (elementary) reduction possible. Conversely, suppose an element of $\mathrm{arr} (\mathcal{B})$ is in normal form. Then it must be some (not \textit{any}) monoidal product and composite of the following ``building blocks'': $$ P_{\sigma} \quad , \quad \mu^{[q]}  \quad , \quad \Delta^{[p]}  \quad , \quad       (\mu \otimes \mu)(\id \otimes P_{A,A} \otimes \id)(\Delta \otimes \Delta) \quad , \quad \eta \otimes \eta  \quad , \quad  \varepsilon  \otimes  \varepsilon$$ (note that any of the three first morphisms includes the identity). By the naturality of the symmetry $P$, any monoidal product and composite of arrows $P_{\sigma}$ and iterated multiplications $\mu^{[q]}$ equals (in the category $\mathcal{B}$)  $ \mu^{[q_1, \ldots , q_m]} \circ P_{\sigma '}  $ for some $\sigma' \in \Sigma_{q_1 + \cdots + q_m}$ and some $m \geq 0$. Similarly, any monoidal product and composite of arrows $P_{\sigma}$ and iterated comultiplications $\mu^{[q]}$ equals  $ P_{\sigma '} \circ \Delta^{[p_1, \ldots , p_n]}   $  for some $n \geq 0$. Furthermore, we can get the three latter arrows by means of the three first arrows, setting $p=0=q$ and using the monoidal product.  Hence we are left to a monoidal product and composite of the three first arrows. But in that composite we cannot find a $\mu^{[q]}$ before a $\Delta^{[p]}$ since it would violate the normal form hypothesis. So we are left to the elements of the given set, as required.
\end{proof}

\begin{proof}[Proof of Theorem \ref{thm normal form}]
By the previous lemma, it suffices to check that the reduction $\rightarrow$ on $\mathrm{arr} (\mathcal{B})$ satisfies the conditions of the Diamond lemma, since there is a canonical full functor $\pi: \mathcal{B} \to \mathsf{B}$ with the property that  any $f_1  \leftarrow f   \rightarrow f_2$ maps to  equalities. Indeed given a morphism $f$ in $\mathsf{B}$, the subset $\pi^{-1}(f) \subset \mathrm{arr}(\mathcal{B})$ is connected as reductions in $\mathcal{B}$ become relations in $\mathsf{B}$. By the Diamond lemma, all elements of $\pi^{-1}(f)$ must have the same (unique) normal form, which by \ref{lem elements in normal form} must be of the form $x= \mu^{[q_1, \ldots , q_m]} \circ P_{\sigma} \circ \Delta^{[p_1, \ldots , p_n]} $. We conclude as $f= \pi(x)$.

Let us then check the conditions of the Diamond lemma. That the  relation $\rightarrow$ satisfies the descending chain condition follows since any element of $\mathrm{arr} (\mathcal{B})$ is a finite composite of finitely many generators, there are finitely many elementary reductions, there are no reduction loops (if $f \rightarrow g$ then no reduction chain starting from $g$ can hit $f$, as $\mathcal{B}$ arises as a free category where  only the (co)associativity relation is modded out) and the identity morphism is in normal form.

To check the diamond condition, let $f_1  \leftarrow f   \rightarrow f_2$. By definition, any element of $\mathrm{arr} (\mathcal{B})$ can be written as a composition of tensor products of the elementary maps $\id$, $ \mu$, $\Delta$, $\eta$, $\varepsilon$ and $P_{1,1}$. Let us argue by induction on the number $n$ of such compositions in $f$, that is, $$ f=(h_{n,1} \otimes \cdots \otimes h_{n, k_n}) \circ \cdots \circ (h_{1,1} \otimes \cdots \otimes h_{1, k_1}), $$ where $h_{i,j}$ is one of the aforementioned elementary maps.

Let us start with the base case $n=2$ (note that for $n=1$ no reduction is possible). The key observation is that  all defining elementary reductions express how a composite of exactly two maps is transformed, namely
\begin{equation}\label{eq:elementary_composit}
\mu (\id \otimes \eta)\ , \ \mu (\eta \otimes \id) \ , \  (\id \otimes \varepsilon)\Delta \ , \  (\varepsilon \otimes \id)\Delta  \ , \  \Delta \mu \ , \  \Delta \eta \ , \  \varepsilon\mu \ , \  \varepsilon \eta.
\end{equation}
Therefore, either $f_1 =f_2$ and the diamond is trivial, or $f$ can be written as a tensor product of maps  with at least two of the arrows \eqref{eq:elementary_composit} as factors, and $f_1$ and $f_2$ are the results of applying the corresponding elementary reduction to $f$ . In this case, the diamond condition also holds by applying the elementary reduction that induces $f \rightarrow f_1$ in the corresponding factor of $f_2$, and conversely for $f_1$, yielding $f_1 \rightarrow g \leftarrow f_2$.

Let us now focus on the case $n=3$. Given $f_1  \leftarrow f   \rightarrow f_2$, if $f_1=f_2$ or the elementary reductions that induce the reductions $f \rightarrow f_i$ appear in different tensor factors, then the same argument as for $n=2$ applies. Otherwise the following hold:
\begin{itemize}
\item $f$ can be written as a tensor product of maps in which one of the factors is given by one of the composites
\begin{equation}\label{eq:comp_3}
 \Delta \mu (\id \otimes \eta) \ ,  \ (\id \otimes \varepsilon)\Delta   \eta   \ ,  \    (\id \otimes \varepsilon) \Delta \mu   \ ,  \    \varepsilon \mu (\id \otimes  \eta), 
\end{equation}
\item the reduction $f   \rightarrow f_1$ is  induced by the elementary reduction given by the composition of the two leftmost maps, e.g. $\Delta \mu$,
\item the reduction $f   \rightarrow f_2$ is  induced by the elementary reduction given by the composition of the two rightmost maps, e.g. $\mu (\id  \otimes \eta)$.
\end{itemize} 
It then suffices to check that if $f$ is any of the four maps \eqref{eq:comp_3}, given a diagram $f_1  \leftarrow f   \rightarrow f_2$ with $f_1$ and $f_2$ as explained above, the diamond condition holds. Since the two latter composites of \eqref{eq:comp_3} are dual to the former two, it is enough to check these. On the one hand we have
\begin{align*}
 \Delta \mu (\id \otimes \eta)  &\rightarrow (\mu \otimes \mu)(\id \otimes P_{1,1} \otimes \id)(\Delta \otimes \Delta)(\id \otimes \eta)\\
 &\rightarrow  (\mu \otimes \mu)(\id \otimes P_{1,1} \otimes \id)(\Delta \otimes \eta \otimes  \eta)\\
 &= (\mu \otimes \mu)(\id \otimes \id  \otimes  \id \otimes \eta)(\id \otimes P_{1,1}  )(\Delta \otimes \eta )\\
 &\rightarrow  (\mu \otimes \id)(\id \otimes P_{1,1} )(\Delta \otimes \eta )\\
 &= (\mu \otimes \id)(\id \otimes \eta \otimes \id )\Delta \\
 &\rightarrow \Delta.
\end{align*}
and on the other hand trivially $ \Delta \mu (\id \otimes \eta)  \rightarrow \Delta.$ Moving on to the second composition of maps in \eqref{eq:comp_3}, we have
$$ (\id \otimes \varepsilon)\Delta   \eta \rightarrow  (\id \otimes \varepsilon)(\eta \otimes \eta) \rightarrow \eta  $$ and on the other hand $ (\id \otimes \varepsilon)\Delta   \eta  \rightarrow \eta$. Any case, we have the diamond condition.

Let us now treat the general case $n >3$. Write $f= g_n \circ g_{n-1} \circ \cdots \circ g_1$, where $g_i=(h_{i,1} \otimes \cdots \otimes h_{i,k_i})$ as above. Given a diagram $f_1  \leftarrow f   \rightarrow f_2$, if the two reductions take place inside the composite $g_{n-1} \circ \cdots \circ g_1$, then we apply the induction hypothesis and conclude. If one of them involves the map $g_n$, then we have two possibilities: if the other reduction takes place in $g_{n-3} \circ \cdots \circ g_1$, then the two reductions cannot interact and the diamond holds trivially. Alternatively, if the other reduction takes place in $g_n \circ g_{n-1} \circ g_{n-2}$, then the diamond follows by the case $n=3$, that we have already shown.
\end{proof}

\begin{corollary}\label{cor normal form algebras}
Every morphism $f: n \to m$ in $\mathsf{A}$ factors in a unique way as $$f= \mu^{[q_1, \ldots , q_m]} \circ P_{\sigma}   $$ for some (unique) integers $q_1, \ldots , q_m \geq 0$ and $\sigma \in \Sigma_n$.
\end{corollary}

\begin{remark}
The existence of the factorisation of  \ref{thm normal form} is shown in  \cite[Lemma 2]{habiro2016} in a more general context (namely for the PROP for commutative Hopf algebras) using a different approach, yet the arguments used there are essentially the same as the ones used in \ref{lem elements in normal form}. We prefer to use an argument involving the diamond lemma to get existence and uniqueness at once. Also note that for the PROP for (commutative) Hopf algebras, uniqueness is not possible, by the antipode axiom.
\end{remark}

\subsection{The main theorem}

Let us write $\N$ for the monoid of non-negative integers. Recall that the forgetful functor $U: \mathsf{CMon} \to \mathsf{Mon}$ from commutative monoids to monoids has a left-adjoint $(-)^{ab} : \mathsf{Mon} \to \mathsf{CMon}$ given by the \textit{abelianisation} of the monoid, that is, the quotient by the smallest congruence relation containing the relations $mn \sim nm$. In particular, we have $$F(x_1, \ldots , x_n )^{ab}\cong \N x_1 \oplus \cdots \oplus \N x_n$$ where $\oplus$ is the biproduct (both product and coproduct) of the category  $\mathsf{CMon}$.

Denote by $\pi_i$ the composite $$F(x_1, \ldots , x_n ) \to \N x_1 \oplus \cdots \oplus \N x_n \to \N x_i$$ where the first map is the natural map to the abelianisation and the second map is the $i$-th projection. Likewise, we will denote as $\pi$ the composite $$F(x_1, \ldots , x_n ) \to \N x_1 \oplus \cdots \oplus \N x_n \to \N$$ where the first map is as before and the second map is the unique map coming from the coproduct such that the composites  $\N x_i \to \N x_1 \oplus \cdots \oplus \N x_n \to \N$ are the identity maps (this is the ``sum'' map). For $w \in  F(x_1, \ldots , x_n )$, it should be clear that $\pi_i (w)$ counts how many times the letter $x_i$ appears in $w$, whereas $\pi(w)$ counts the total number of letters $w$ contains.


\begin{definition}
Let $\widehat{\mathsf{fgFMon}}$ be the following category: its objects are the free monoids $F([n]) = F(x_1, \ldots , x_n)$, $n \geq 0$. An arrow is a pair $(f, \bm{\sigma})$, where $$f: F(x_1, \ldots , x_n)\to F(y_1, \ldots , y_m)$$ is a monoid map and $\bm{\sigma} = (\sigma_1, \ldots , \sigma_m)$ is a tuple of $m$ permutations $\sigma_i \in \Sigma_{k_i}$, where $k_i := \pi_i (f(x_1 \cdots x_n))$. 
The identity of an object $ F(x_1, \ldots , x_n)$ is given by the identity monoid homomorphism and $n$ copies of the unique element $\id_1 \in \Sigma_1$. Moreover, this category is symmetric monoidal (hence a PROP) where the monoidal product is given by the free product $*$ of monoids on objects and  $$(f,\bm{\sigma}) \otimes (g, \bm{\tau} ) := (f * g, \bm{\sigma \tau})$$ on morphisms. In the above expression, $\bm{\sigma \tau}$ stands for the concatenation of the two tuples of permutations. The unit is the trivial monoid $\bm{1} = F([0])$.

We still have to describe the composite law in $\widehat{\mathsf{fgFMon}}$. Setting up a formula for the composite is a highly non-trivial task and requires several constructions that we describe in the next subsection. Yet we would like to place the expression here for future reference. Given arrows $(f, \bm{\sigma}): F(x_1, \ldots x_n) \to F(a_1,  \ldots , a_m)$ and $(g, \bm{\tau}):  F(a_1,  \ldots , a_m) \to F(s_1, \ldots, s_\ell)$, the composite  $(g, \bm{\tau})\circ (f, \bm{\sigma}) =: (g \circ f, \bm{\tau} \circ \bm{\sigma})$ will consist of $\ell$ permutations $$\bm{\tau} \circ \bm{\sigma} := ((\bm{\tau} \circ \bm{\sigma})_1, \ldots , (\bm{\tau} \circ \bm{\sigma})_\ell)$$ where $ (\bm{\tau} \circ \bm{\sigma})_i \in \Sigma_{\pi_i (gf (x_1 \cdots x_n))}$. These permutations are uniquely determined by the following formula:
\begin{align}\label{eq composite law fgFMonhat}
\begin{split}
(\bm{\tau} \circ \bm{\sigma})_1 \otimes \cdots \otimes (\bm{\tau} \circ \bm{\sigma})_\ell &= \xi_{(g\circ f)(x_1 \cdots x_n)}\circ \rho_{f, \mu^{[m]} \circ g} \circ (\gamma_{(k_1, p_1)} \otimes \cdots \otimes  \gamma_{(k_m, p_m)})\\
 &\circ (\sigma_1^{\otimes p_1} \otimes \cdots \otimes \sigma_m^{\otimes p_m}) \circ \langle  \Psi (g(a_1 \cdots a_m), \bm{\tau}) \rangle_{k_1^{\times p_1}, \ldots k_m^{\times p_m}}
\end{split}
\end{align}
where $p_i:= \pi (g (a_i))$ and $k_i := \pi_i (f (x_1 \cdots x_n))$ for $1 \leq i \leq m$. The symbols used here are explained in \ref{constr Phi}, \ref{constr xi gamma} and \ref{constr Psi}.
\end{definition}

Although strange, the previous composite law yields the following

\begin{theorem}\label{thm fundamental theorem}
The categories $\widehat{\mathsf{fgFMon}}$ and $ \mathsf{B}$ are monoidally equivalent. In other words, they are equivalent PROPs.
\end{theorem}

The proof of this theorem consists of a long construction and will be carried out in next subsection. Part of this construction will consist of establishing  the expression \eqref{eq composite law fgFMonhat}.

\begin{remark}
Before presenting the proof, we would like to guide the reader through an informal argument supporting the theorem presented above. Starting from the base that there is an equivalence $\mathsf{ComB} \simeq \mathsf{fgFMon}$ as explained in \ref{ex ComB = fgFMon}, the permutations that we added to the monoid maps will fix the order of the multiplication. More concretely, $(\mu, \id) \circ (P, \bm{\id}) = (\mu, (12))$, compare with \ref{ex motivating 1}. Yet these permutations do not add extra data to the rest of structure maps: for the iterated comultiplication $\Delta^{[m]}: F(x) \to F(y_1, \ldots , y_m)$, $x \mapsto y_1 \cdots y_m$, we have that $\pi_i (y_1 \cdots y_m) =1$, so the $m$ permutations that accompany  $\Delta^{[m]}$ must be the trivial ones, $\id_1 \in \Sigma_1$. For the unit $\eta: \bm{1} \to F(y)$ and the counit $\varepsilon : F(x) \to 1$ this is even more vacuous, since there cannot be any permutation attached. Thence the permutations only add extra data for the multiplication map, which is precisely what we want to refine.
\end{remark}

\subsection{Proof of Theorem \ref{thm fundamental theorem}}
Our strategy consists of constructing, for any symmetric monoidal category $\C$, a bijection  $$\mathsf{Alg}_{\widehat{\mathsf{fgFMon}}} (\C) \toiso  \mathsf{Alg}_{\mathsf{B}} (\C)$$  natural on $\C$ and appeal to \ref{lemma yoneda}. The hardest part of the proof will be to check (or rather, to define) that the cumbersome composite law \eqref{eq composite law fgFMonhat} models the order of multiplication in a bialgebra. The first step in the construction is the following 


\begin{lemma}\label{lem F(x) bialgebra object}
The object  $ F([1])=F(x) \in \widehat{\mathsf{fgFMon}}$ is a bialgebra object. It has structure maps  $\mu, \eta, \Delta, \varepsilon$ whose underlying monoid maps are the same as in   \ref{ex ComB = fgFMon} and whose attached permutations are the identities.
\end{lemma}

The proof of the lemma requires to handle the composition law of $\widehat{\mathsf{fgFMon}}$, that we agreed we will describe later. Thus we defer this proof until page \pageref{proof lem F(x) bialg}.

Let us address the construction of the desired bijection. By \ref{lem F(x) bialgebra object}, there exists a unique strong monoidal functor $$T: \mathsf{B} \to \widehat{\mathsf{fgFMon}}$$ that maps the generating object of $\mathsf{B}$ to the object $F(x)$. Precomposition along this functor defines a map 
\begin{equation}\label{eq bijection alg_B = alg_fgFMonhat}
\mathsf{Alg}_{\widehat{\mathsf{fgFMon}}} (\C) \to  \mathsf{Alg}_{\mathsf{B}} (\C)
\end{equation}
natural on $\C$. 

We aim to define an inverse for the previous map \eqref{eq bijection alg_B = alg_fgFMonhat}. For that purpose a key point will be to understand how we can get large permutations from smaller ones ``patched'' according to a certain element of a free monoid.\medskip



Recall that a \textit{partition} in $m\geq 0$ blocks on a set $S$ is a collection of $m$ pairwise disjoint subsets $S_1, \ldots , S_m  \subseteq S$ (some of them possibly empty) such that $\cup_i S_i =S$. If $S$ is ordered, so are the subsets $S_i$. We denote by $\mathrm{Part}_m (S)$ the set of partitions in $m$ blocks on $S$.

\begin{construction}\label{constr Phi}
Let $m \geq 0$. Given a word $w \in F(y_1, \ldots , y_m) $, we will construct a partition of the set $[\pi (w)]$ given by $m$ subsets $S_i$, each of them with $\pi_i(w)$ elements. Write $w$  as a product of the elements $y_i$ with no powers (this is a product of $\pi (w)$ letters). Then $S_i$ is the ordered subset containing the positions, from left to right, where the letter $y_i$ appears in $w$.  Denote by $\Phi (w)$ the partition defined this way.
\end{construction}

\begin{lemma}
Let $m \geq 0$. The previous construction defines a bijection
\begin{equation}
\Phi : F(y_1, \ldots , y_m) \toiso  \coprod_{N \geq 0}  \mathrm{Part}_m ([N]).
\end{equation}
\end{lemma}
\begin{proof}
The obvious inverse is  defined as the map taking a partition $S_1,  \ldots , S_m$ of $[N]$ to the word $w$ product of $N$ letters $y_i$ according to the order specified by the subsets $S_i$. This establishes the bijection.
\end{proof}

\begin{construction}\label{constr xi gamma}
Given a word $w \in F(y_1, \ldots , y_m)$ and $\Phi(w)=(S_1, \ldots , S_m)$, let $\xi_w \in  \Sigma_{\pi(w)}$ be the unique permutation mapping $S_i$ to $$\left\lbrace \sum_{j=1}^{i-1} \pi_j(w) + 1 < \cdots < \sum_{j=1}^{i} \pi_j(w) \right\rbrace  $$  in a order-preserving way. We can think of this permutation  as ``turning'' the word $a_{1}^{\pi_1(w)} \cdots a_{m}^{\pi_m(w)}$ into $w$. A schematic of this permutation is depicted in Figure \ref{fig omega}. A special instance of this construction is the case where $w_0= (y_1 \cdots y_m) \overset{p}{\cdots}(y_1 \cdots y_m)$. In this case we denote $$\gamma_{(m,p)}:= \xi_{w_0}^{-1} \in \Sigma_{mp}.$$ Explicitly, $$\gamma_{(m,p)} (lm+r) := (r-1)p+(l+1)$$ for $0 \leq l \leq p-1$ and $1 \leq r \leq m$. It is straightforward to check that $\gamma_{(m,1)} = \id_m $ and $\gamma_{(1,p)} = \id_p$.

Note that if
\begin{equation}\label{eq w k_i,j}
w= y_1^{k_{1,1}} y_2^{k_{2,1}} \cdots y_m^{k_{m,1}} y_1^{k_{1,2}}y_2^{k_{2,2}} \cdots y_m^{k_{m,p}} \in F(y_1, \ldots , y_m)
\end{equation}
for some $k_{i,j} \geq 0$, then we have  $$ \xi_{w}= \langle \xi_{w_0} \rangle_{k_{1,1}, k_{1,2}, \ldots , k_{p,m}} \in \Sigma_{\pi(w)}.$$

\begin{figure}[h]
\begin{tikzpicture}
\draw (0, 0) node {$\overbrace{\underbrace{k_{1,1}} \  \underbrace{k_{1,2}} \quad \cdots \quad \underbrace{k_{1,p}}}^p \quad \overbrace{\underbrace{k_{2,1}} \quad \cdots \quad \underbrace{k_{2,p}}}^p \qquad \cdots \qquad \overbrace{\underbrace{k_{m,1}} \quad \cdots \quad \underbrace{k_{m,p}}}^p $};
\draw (0, -3) node {$\underbrace{\overbrace{k_{1,1}}  \  \overbrace{k_{2,1}} \quad \cdots \quad \overbrace{k_{m,1}}}_m \quad \underbrace{\overbrace{k_{1,2}} \quad \cdots \quad \overbrace{k_{m,2}}}_m  \qquad \cdots \qquad \underbrace{\overbrace{k_{1,p}} \quad \cdots \quad \overbrace{k_{m,p}}}_m $};
\draw [->] (-5,-0.7) -- (-5,-2.3);
\draw [->] (-1.4,-0.7) to [bend right] (-4.2,-2.3);
\draw [->] (3,-0.7) to [bend right] (-2.5,-2.3);
\draw [->] (-4.2,-0.7) to [bend right] (-1.5,-2.3);
\draw [->] (-2.5,-0.7) to [bend right] (3,-2.3);
\draw [->] (5,-0.7) -- (5,-2.3);
\end{tikzpicture}
\caption{Schematic of the permutation $\xi_w \in \Sigma_{\pi(w)}$. Given an arbitrary word $w \in F(y_1, \ldots , y_m)$ of the form \eqref{eq w k_i,j}, if we divide $\{ 1 <  \cdots < \pi(w) \}$ in $mp$ blocks of sizes $k_{i,j}$, then $\xi_w$ shuffles them placing the blocks $k_{1,1}, \ldots , k_{m,1}$ in the first position, after that the blocks $k_{1,2}, \ldots , k_{m,2}$, etc. If $k_{i,j} =1$ for all $i,j$, then this illustration represents $\gamma_{(m,p)}^{-1}$.}
\label{fig omega}
\end{figure}
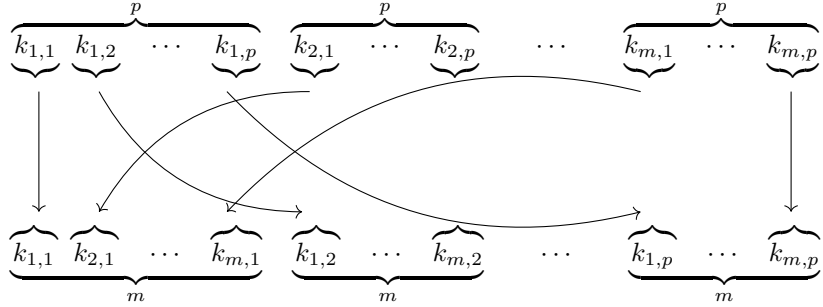
\end{construction}

If $f: F(x) \to F(y_1, \ldots , y_m)$ with $f(x)=w$ and $g: F(y_1, \ldots , y_m) \to F(s)$, with $f(y_i)= s^{p_i}$, $p_i \geq 0$, we define
\begin{equation}
\rho_{f,g}  := (\langle \xi_w  \rangle_{p_1^{\times k_1}, \ldots ,  p_m^{\times k_m}})^{-1}
\end{equation}
where $k_i:= \pi_i (w)$ and $p_1^{\times k_1}$ stands for a sequence $p_i , \ldots , p_i$ with $k_i$ times $p_i$.

\begin{construction}\label{constr Psi}
Let $k_1, \ldots k_m \geq 0$ and let $n= \sum_i k_i$. Write $$\mathcal{W}(k_1, \ldots , k_m):= \{ w \in F(y_1, \ldots, y_m) : \pi_i (w) = k_i , \ 1 \leq i \leq m  \}.$$ 
We define a map 
$$\Psi_{k_1, \ldots , k_m} = \Psi : \mathcal{W}(k_1, \ldots , k_m) \times \Sigma_{k_1} \times \cdots \times \Sigma_{k_m} 	\to \Sigma_n$$ as follows: 
%
%
%
%
%
$$\Psi (w, \sigma_1, \ldots , \sigma_m) := \xi_w^{-1} \circ  (\sigma_1 \otimes \cdots \otimes \sigma_m)   .$$
\end{construction}

\begin{lemma}\label{lem utilisimo}
The previous map $$ \Psi: \mathcal{W}(k_1, \ldots , k_m) \times \Sigma_{k_1} \times \cdots \times \Sigma_{k_m} \toiso  \Sigma_n $$ is a bijection. Furthermore, it is ``monoidal'' in the sense that the diagram
$$
\begin{tikzcd}
\mathcal{W}(k_1, \ldots , k_m) \times \prod_{i=1}^m \Sigma_{k_i} \times \mathcal{W}(k_1', \ldots , k'_r) \times \prod_{i=1}^r \Sigma_{k'_i} \rar{\Psi_{k_1, \ldots , k_m} \times \Psi_{k'_1, \ldots , k'_r}} \dar &[3.8em]  \Sigma_n  \times  \Sigma_{n'} \dar{\otimes} \\
\mathcal{W}(k_1, \ldots , k_m, k'_1, \ldots , k'_r) \times \prod_{i=1}^m \Sigma_{k_i} \times \prod_{i=1}^r \Sigma_{k'_i} \rar{\Psi_{k_1, \ldots , k_m, k'_1, \ldots , k'_r}} & \Sigma_{n +n'}
\end{tikzcd}
$$
commutes, where $n'= \sum_i k'_i$ and the left-hand arrow concatenates the words and acts as the identity on the symmetric groups.
\end{lemma}
\begin{proof}
Observe that the cardinality of $\mathcal{W}(k_1, \ldots , k_m)$ is the multinomial coefficient $\binom{n}{ k_1, \ldots , k_m}$, which means that the source and target of $\Psi$ have the same cardinality since $$  \binom{n}{ k_1, \ldots , k_m} = \frac{n!}{k_1 ! \ldots k_m !}.$$ Hence it suffices  to check that $\Psi$ is surjective: given $\alpha \in \Sigma_n$, turn the formal word $z_{1} \cdots z_{n}$ into a word $w \in F(y_1, \ldots, y_m)$ by the following rule: $$z_{i} := \begin{cases}
y_1, & 1 \leq \alpha^{-1} (i) \leq k_1\\
y_2, & k_1+1 \leq \alpha^{-1}  (i) \leq k_2\\
\vdots \\
y_m, & \sum_{j=1}^{m-1} k_j+1 \leq \alpha^{-1} (i) \leq \sum_{j=1}^{m} k_j
\end{cases}$$ It is clear that $w \in \mathcal{W}(k_1, \ldots , k_m)$. By construction,  $   \xi_w \circ \alpha$ is of the form $\sigma_1 \otimes \cdots \otimes \sigma_m$ for some (unique) permutations $\sigma_i \in \Sigma_{k_i}$, what establishes the bijection.

The commutativity of the diagram follows from the observation that the block product of permutations is compatible with the group law of the symmetric group\footnote{A more concise way to restate this is by saying that if we view $\Sigma_n$ as a groupoid in the usual way, then the block product $\otimes$ makes $\Sigma_n$ a monoidal category.} in the sense that $$(\sigma \otimes \sigma') \circ (\tau \otimes \tau') = (\sigma \circ \tau) \otimes (\sigma' \circ \tau')   .$$ Then we have
\begin{align*}
\Psi_{k_1, \ldots , k_m} (w, \sigma_1, \ldots , \sigma_m) &\otimes \Psi_{k'_1, \ldots , k'_r} (w', \sigma'_1, \ldots , \sigma'_m) \\ &= [\xi_w^{-1} \circ  (\sigma_1 \otimes \cdots \otimes \sigma_m)] \otimes [\xi_{w'}^{-1} \circ  (\sigma'_1 \otimes \cdots \otimes \sigma'_r)]\\
&= (\xi_w^{-1} \otimes \xi_{w'}^{-1}) \circ (\sigma_1 \otimes \cdots \otimes \sigma_m \otimes \sigma'_1 \otimes \cdots \otimes \sigma'_r)\\
&= \xi_{ww'}^{-1}  \circ (\sigma_1 \otimes \cdots \otimes \sigma_m \otimes \sigma'_1 \otimes \cdots \otimes \sigma'_r)\\
&= \Psi_{k_1, \ldots , k_m, k'_1, \ldots , k'_r} (ww', \sigma_1, \ldots , \sigma_m, \sigma'_1, \ldots , \sigma'_m)
\end{align*}
which concludes.
\end{proof}

Now we are ready to define an inverse for \eqref{eq bijection alg_B = alg_fgFMonhat}.
\begin{definition}
 Given a bialgebra $A \in \C$, define a functor $$\widehat{A}: \widehat{\mathsf{fgFMon}} \to \C$$ as follows: on objects, $\widehat{A}(F(x_1, \ldots , x_n)):= A^{\otimes n}$. On morphisms, given an arrow $(f, \bm{\sigma}): F(x_1, \ldots , x_n) \to F(y_1, \ldots , y_m)$, let
$$ p_i := \pi(f(x_i)) \quad , \quad q_j := \pi_j (f(x_1 \cdots x_n)) \quad , \quad \sigma := \Psi (f(x_1 \cdots x_n), \bm{\sigma}) $$ for $1 \leq i \leq n$ and $1 \leq j \leq m$ and define
\begin{equation}\label{eq def Ahat}
\widehat{A} (f, \bm{\sigma}):= \mu^{[q_1, \ldots , q_m]} \circ P_{\sigma} \circ \Delta^{[p_1, \ldots , p_n]} . 
\end{equation}
\end{definition}

Let us illustrate the above assignment.

\begin{example}
Let $(f, \bm{\sigma}): F(x,y) \to F(a,b)$ be the map in $\widehat{\mathsf{fgFMon}}$ determined by $$  f(x):= a^2 b \quad , \quad f(y):= abab \quad , \quad  \sigma_1 := (4321) \quad , \quad  \sigma_2 := (13),$$ where $\bm{\sigma} = (\sigma_1, \sigma_2)$. Note that $\sigma_1 \in \Sigma_4$ as $4= \pi_1 (f(xy)) = \pi_1 (a^2babab)$ and similarly $\sigma_2 \in \Sigma_3$. We easily see that $p_1 = 3$, $p_2= 4$, $q_1 =4$, $q_2 = 3$. Then
\begin{align*}
\sigma &= \Psi (a^2babab, (4321), (13)) = \xi_{ab^2abab}^{-1} \circ [(4321) \otimes (13)] \\ &= (3456)(4321)(57)= (165732) \in \Sigma_7.
\end{align*}
Alternatively, this can also be computed as follows: consider the ordered sets $A:= \{ z_1, z_2, z_3, z_4 \}$ and $B:= \{ z_5, z_6, z_7 \}$, and let $\sigma_1$ act on $A$ and $\sigma_2$ on $B$, which yields $\sigma_1 A = \{z_4, z_1, z_2, z_3  \}$ and $  \sigma_2 B = \{ z_7, z_6, z_5  \}$. Next consider the ordered union $\sigma_1 A \amalg \sigma_2 B = \{z_4, z_1, z_2, z_3 , z_7, z_6, z_5 \}$. Now for $1 \leq i \leq 4$, we replace the letter $z_i$ by $z_{o(i)}$ where $o(i)$ is the position of the $i$-th letter $a$ in the pattern $a a b abab$. Similarly for $1 \leq i \leq 3$ replace  $z_{i+4}$ by $z_{o'(i)}$ where $o'(i)$ is the position of the $i$-th letter $b$ in the pattern, according to the order from left to right. This yields the pattern $z_6 z_1 z_2 z_4 z_7 z_5 z_3$. Then $\sigma$ is the permutation such that the previous pattern equals $z_{\sigma (1)} \cdots z_{\sigma(7)}$, that is, $\sigma = (165732) $, as above.

Following \ref{notation} we conclude that $\widehat{A} (f, \bm{\sigma})$ is the map $$x \otimes y \mapsto y_{(2)} x_{(1)} x_{(2)} x_{(4)} \otimes y_{(3)} y_{(1)} x_{(3)}.$$
\end{example}

Let us continue with the proof of \ref{thm fundamental theorem}. The key point will be to set the expression \eqref{eq composite law fgFMonhat}, which has to model all changes of indices that take place when applying a combination of multiplication and comultiplication maps to an element of a bialgebra. In other words, we have to set \eqref{eq composite law fgFMonhat} so that  $\widehat{A}$ actually defines a functor.


\begin{proposition}
Given a bialgebra $A \in \C$, we have that
\begin{equation}\label{eq what is the composite}
 \widehat{A}  ( (g, \bm{\tau}) \circ (f, \bm{\sigma})  )  = \widehat{A}(g, \bm{\tau})  \circ \widehat{A}(f, \bm{\sigma}) ,
\end{equation}
so that  $\widehat{A}: \widehat{\mathsf{fgFMon}} \to \C$ is indeed a functor.
\end{proposition}
\begin{proof}
Instead of verifying \eqref{eq what is the composite}  using \eqref{eq composite law fgFMonhat}, the proof will construct the composite law \eqref{eq composite law fgFMonhat} of $\widehat{\mathsf{fgFMon}}$ by \textit{imposing} the equality \eqref{eq what is the composite} and the associativity of the composite. We will subdivide the proof in several steps, increasing the complexity of the maps at hand.



\noindent \textit{Step 1.} We start by composing maps between the objects $F([1])$. Let $(f, \sigma): F(x) \to F(a)$, $x \mapsto a^k$ and $(g, \tau): F(a) \to F(s)$, $a \mapsto s^p$ for integers $k,p \geq 0$. If $\sigma= \id_k \in \Sigma_k $ and $\tau = \id_p \in \Sigma_p$, then their images under $\widehat{A}$ correspond to the maps $x \mapsto x_{(1)} \cdots x_{(k)}$ and $a \mapsto a_{(1)} \cdots a_{(p)}$. Their composite  $\widehat{A}(g, \id_p)  \circ \widehat{A}(f, \id_k)$ is
\begin{equation}\label{eq base case}
x \mapsto (x_{(1)} \cdots x_{(k)})_{(1)} \cdots ( x_{(1)} \cdots x_{(k)})_{(p)} = \prod_{i=1}^p x_{(i)} x_{(p+i)} \cdots x_{((n-1)p+i)}.
\end{equation}
This implies that the permutation associated to the  composite is given by the element ${\gamma_{(k,p)} \in  \Sigma_{kp}}$, which was already defined in \ref{constr xi gamma}. If $\sigma, \tau$ are not the identity permutations, then the left-hand side of the above expression becomes $$x \mapsto (x_{(\sigma 1)} \cdots x_{(\sigma k)})_{(\tau 1)} \cdots ( x_{(\sigma 1)} \cdots x_{(\sigma k)})_{(\tau p)} .$$ This says that $\tau$ acts on \eqref{eq base case} permuting blocks of size $k$, whereas $\sigma$ acts inside each of the $p$ blocks. All in all, the permutation associated to the composite $\widehat{A}  ( (g, \tau) \circ (f,\sigma)  )$ must be
\begin{equation}
 \gamma_{(k,p)} \circ  (\sigma \otimes \overset{p}{\cdots} \otimes \sigma) \circ  \langle \tau \rangle_{k,  \overset{p}{\cdots} , k} \in \Sigma_{kp}.
\end{equation}

\noindent \textit{Step 2.} We want to study now the case where $(f, \bm{\sigma}): F(x) \to F(a_1,  \ldots , a_m)$, $x \mapsto w$,  $\bm{\sigma} = (\sigma_1, \ldots, \sigma_m)$ and $(g, \tau):  F(a_1,  \ldots , a_m) \to F(s)$, $a_i \mapsto  s^{p_i}$, $p_i \geq 0$. Let us start from the base case where $\sigma_i=  \id$ , $\tau=\id$ and $w = a_1^{k_1} \cdots a_m^{k_m}$. Here the situation is almost identical to the first part of Step 1:  $\widehat{A}  (f,\bm{\sigma}) $ corresponds to the map $x \mapsto x_{(1)} \cdots x_{(k_1)} \otimes  x_{(k_1+1)} \cdots x_{(k_2)} \otimes \cdots$ whereas $\widehat{A}  (g,\tau)$ corresponds to $a_1 \otimes \cdots \otimes a_m \mapsto (a_1)_{(1)} \cdots  (a_1)_{(p_1)} \cdots  (a_m)_{(1)} \cdots  (a_m)_{(p_m)}$, so it follows reasoning with each of the ``monoidal blocks'' of $\widehat{A}  (f,\bm{\sigma}) $ that the composite gets as permutation $$\gamma_{(k_1,p_1)}\otimes  \cdots \otimes  \gamma_{(k_m,p_m)}  \in \Sigma_{N}$$ where  $N= k_1p_1 + \cdots + k_mp_m$.  Now suppose that $f(x)=w  \in F(a_1, \ldots , a_m)$ is a general element.  We want to show that if $ \rho_{f,g}$ is defined from $\Phi (w)$ as above, the permutation associated to the composite $\widehat{A}  ( (g, \tau) \circ (f,\bm{\sigma})  )$ is $$ \rho_{f,g}  \circ ( \gamma_{(k_1,p_1)}\otimes  \cdots \otimes  \gamma_{(k_m,p_m)} ) \in \Sigma_{N}.$$ For observe that the permutation $\xi_w$ is designed to keep track of the ``reordering'' of the word $a_{1}^{\pi_1(w)} \cdots a_{m}^{\pi_m(w)}$ turning into $w$. When passing to $ \langle \xi_w  \rangle_{p_1^{\times k_1}, \ldots ,  p_m^{\times k_m}} = \rho_{f,g}^{-1}$, we extend this permutation to a block permutation according to the $p_i$'s. Taking the inverse is needed (just as with $\gamma_{(m,p)}$) to match the indices and not the positions. Lastly, if $\bm{\sigma}$ and $\tau$ are not trivial, then the situation is again similar to the one in Step 1: we have that each of the $ \sigma_i$ acts inside $p_i$ blocks (of size  $m_i$), whereas $\tau$ acts permuting these blocks. An illustration of this is in Figure \ref{fig permutations}.

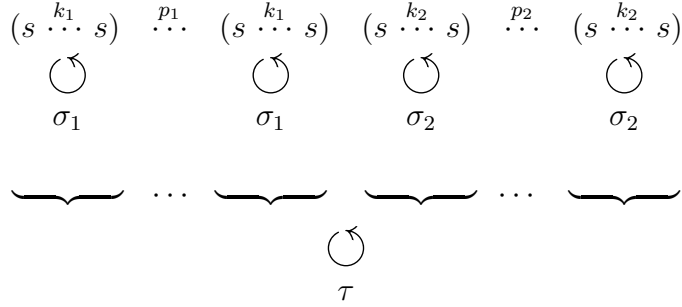
\begin{figure}[h]
\begin{tikzpicture}
\Large \draw (0, 0) node {$(s \  \overset{k_1}{\cdots} \  s) \quad \overset{p_1}{\cdots}  \quad (s \  \overset{k_1}{\cdots} \  s) \quad(s \  \overset{k_2}{\cdots} \  s) \quad \overset{p_2}{\cdots}  \quad(s \  \overset{k_2}{\cdots} \  s)$};
\draw (-3.7,-0.7) node {\Huge $ \circlearrowleft$};
\draw (-3.7,-1.3) node { $\sigma_1$};
\draw (-1,-0.7) node {\Huge $ \circlearrowleft$};
\draw (-1,-1.3) node { $\sigma_1$};
\draw (1,-0.7) node {\Huge $ \circlearrowleft$};
\draw (1,-1.3) node { $\sigma_2$};
\draw (3.7,-0.7) node {\Huge $ \circlearrowleft$};
\draw (3.7,-1.3) node { $\sigma_2$};
\draw (-3.7, -2) node {$\underbrace{\phantom{(s \  \overset{k_1}{\cdots} \  s)}}$};
\draw (-1, -2) node {$\underbrace{\phantom{(s \  \overset{k_1}{\cdots} \  s)}}$};
\draw (-2.3, -2.3) node {$\cdots$};
\draw (2.3, -2.3) node {$\cdots$};
\draw (1, -2) node {$\underbrace{\phantom{(s \  \overset{k_1}{\cdots} \  s)}}$};
\draw (3.7, -2) node {$\underbrace{\phantom{(s \  \overset{k_1}{\cdots} \  s)}}$};
\draw (0, -3) node {\Huge $ \circlearrowleft$};
\draw (0,-3.6) node { $\tau$};
\end{tikzpicture}
\caption{Illustration of  the action of $\bm{\sigma}$ and $\tau$, for $m=2$.}
\label{fig permutations}
\end{figure}

All in all, we conclude the permutation associated to the composite $\widehat{A}  ( (g, \tau) \circ (f,\bm{\sigma})  )$ must be
\begin{equation}\label{eq composite law Step 2}
  \rho_{f,g}  \circ ( \gamma_{(k_1,p_1)}\otimes  \cdots \otimes  \gamma_{(k_m,p_m)} ) \circ  (\sigma_1^{ \otimes p_1}  \otimes \cdots \otimes\sigma_m^{ \otimes p_m})  \circ  \langle \tau \rangle_{k_1^{\times p_1}, \ldots k_m^{\times p_m}} \in \Sigma_{N}
\end{equation}
where $k_i =  \pi_i (w)$ and $N= k_1p_1 + \cdots + k_mp_m.$\medskip

\noindent \textit{Step 3.} Next we consider the case where  $(f, \bm{\sigma}): F(x_1, \ldots , x_n) \to F(a_1,  \ldots , a_m)$, $x_i \mapsto w_i$,  $\bm{\sigma} = (\sigma_1, \ldots, \sigma_m)$  and $(g, \tau):  F(a_1,  \ldots , a_m) \to F(s)$, $a_i \mapsto  s^{p_i}$, $p_i \geq 0$ as before. In this case, the claim is that the formula \eqref{eq composite law Step 2} obtained in Step 2 still holds, taking $w := f(x_1 \cdots x_n)= w_1 \cdots w_n$. Vaguely speaking, the reason is that  in this case $(f, \bm{\sigma})$ models a map $A^{\otimes n} \to A^{\otimes m}$, which is determined by the image of elements $x_1 \otimes \cdots \otimes x_n$. 

More precisely, let us show that if we impose \eqref{eq what is the composite}, then the equation 
\begin{equation}\label{eq imposed associativity 1}
(f, \bm{\sigma}) \circ (\Delta^{[n]}, \bm{\id}) = (f \circ \Delta^{[n]},\bm{\sigma} )
\end{equation}
must hold. We have that 
\begin{align*}
\widehat{A}(f, \bm{\sigma}) \circ \widehat{A}(\Delta^{[n]}, \bm{\id}) &=  (\mu^{[q_1, \ldots , q_m]} \circ P_{\sigma} \circ \Delta^{[p_1, \ldots , p_n]} ) \circ \Delta^{[n]} \\
&= \mu^{[q_1, \ldots , q_m]} \circ P_{\sigma} \circ \Delta^{[p_1+ \cdots + p_n]} .
\end{align*}
Since we want the equation
$$\widehat{A}(f, \bm{\sigma}) \circ \widehat{A}(\Delta^{[n]}, \bm{\id}) = \widehat{A}(f \circ \Delta^{[n]},\bm{\sigma} \circ \bm{\id} )$$
to hold and $\sigma = \Psi (f(x_1 \cdots x_n), \bm{\sigma}) = \Psi ((f \circ \Delta^{[n]})(x), \bm{\sigma}) $, then $\bm{\sigma} \circ \bm{\id}$ must be $\bm{\sigma}$. By associativity of the composite, we must have
\begin{align*}
(g \circ (f \circ \Delta^{[n]}), \tau \circ \bm{\sigma}) &= (g, \tau) \circ [ (f, \bm{\sigma}) \circ (\Delta^{[n]}, \bm{\id})  ]\\
&= [(g, \tau) \circ  (f, \bm{\sigma})] \circ (\Delta^{[n]}, \bm{\id})\\
&= (g \circ f , \tau \circ \bm{\sigma})\circ (\Delta^{[n]}, \bm{\id})\\
&= ((g \circ f) \circ \Delta^{[n]}, \tau \circ \bm{\sigma})
\end{align*}
where we have applied \eqref{eq imposed associativity 1} in the first and third equality. Since the left-hand side of the equation is computable in terms of the formula given in Step 2,  the claim follows.\medskip

\noindent \textit{Step 4.} Let us finally deal with the general case where  $(f, \bm{\sigma}): F(x_1, \ldots x_n) \to F(a_1,  \ldots , a_m)$ and $(g, \bm{\tau}):  F(a_1,  \ldots , a_m) \to F(s_1, \ldots, s_\ell)$, $a_i \mapsto  w'_i$ and $\bm{\tau} =( \tau_1, \ldots , \tau_\ell )$. By Step 3, it is enough to consider $(f, \bm{\sigma}): F(x) \to F(a_1,  \ldots , a_m)$, $x \mapsto w$,  $\bm{\sigma} = (\sigma_1, \ldots, \sigma_m)$. As explained before, the composite $(g, \bm{\tau})\circ (f, \bm{\sigma}) = (gf, \bm{\tau} \circ \bm{\sigma})$ will have $\ell$ associated  permutations that are denoted as $(\bm{\tau} \circ \bm{\sigma})_1, \ldots , (\bm{\tau} \circ \bm{\sigma})_\ell$.\medskip

The strategy we follow is ``dual'' to the one of Step 3, where we read off the composite with $(\Delta^{[n]}, \bm{\id})$. This time we want to show that if $(f,\bm{\sigma}): F(x_1, \ldots , x_n) \to F(a_1, \ldots , a_m)$ then the equation
\begin{equation}\label{eq imposed associativity 2}
(\mu^{[m]}, \bm{\id}) \circ (f,\bm{\sigma}) = (\mu^{[m]} \circ f, \Psi (f(x_1 \cdots x_n), \bm{\sigma}))
\end{equation}
must hold. Indeed for any bialgebra $A \in \C$ we have 
\begin{align*}
\widehat{A}(\mu^{[m]}, \bm{\id}) \circ \widehat{A}(f, \bm{\sigma}) &= \mu^{[m]} \circ (\mu^{[q_1, \ldots , q_m]} \circ P_{\sigma} \circ \Delta^{[p_1, \ldots , p_n]} )  \\
&= \mu^{[q_1 +  \cdots + q_m]} \circ P_{\sigma} \circ \Delta^{[p_1, \ldots , p_n]}
\end{align*}
where $\sigma = \Psi (f(x_1 \cdots x_n), \bm{\sigma})$. Since we want the equation
$$\widehat{A}(\mu^{[m]}, \bm{\id}) \circ \widehat{A}(f, \bm{\sigma}) = \widehat{A}(\mu^{[m]} \circ f ,\bm{\id}\circ \bm{\sigma}  )$$
to hold, then $\bm{\id}\circ \bm{\sigma}$ must be\footnote{Alternatively, we could use \eqref{eq composite law Step 2}, which yields the same result.} $\sigma = \Psi (f(x_1 \cdots x_n), \bm{\sigma}) = \xi_{f(x_1 \cdots x_n)}^{-1} \circ (\sigma_1 \otimes \cdots \otimes \sigma_m)$. In particular, this says that $\bm{\sigma}$ is determined by $\sigma$, since
\begin{equation}\label{eq sigma determines boldsigma}
\sigma_1 \otimes \cdots \otimes \sigma_m = \xi_{f(x_1 \cdots x_n)} \circ \sigma
\end{equation}
and the map $$\Sigma_{k_1} \times \cdots \times \Sigma_{k_m} \hooklongrightarrow \Sigma_{k_1 + \cdots + k_m} \quad , \quad (\sigma_1, \ldots, \sigma_m ) \mapsto \sigma_1 \otimes \cdots \otimes \sigma_m $$ is injective. Therefore, for the composite $ (g, \bm{\tau})\circ (f, \bm{\sigma})$, using the associativity of the composite, we must have
\begin{align*}
((\mu^{[m]} \circ g) \circ f , \Psi (g(a_1 \cdots a_m), \bm{\tau}) \circ \bm{\sigma}) &=   [ (\mu^{[m]}, \bm{\id}) \circ (g, \bm{\tau})] \circ (f, \bm{\sigma}) \\
&= (\mu^{[m]}, \bm{\id}) \circ [ (g, \bm{\tau}) \circ (f, \bm{\sigma})]\\
&= (\mu^{[m]}, \bm{\id}) \circ (g \circ f, \bm{\tau} \circ \bm{\sigma}).
\end{align*}

The first term is computable in terms of Step 2: if $p_i= \pi (g (a_i))$ and $k_i = \pi_i (f (x))$ for $1 \leq i \leq m$, then
\begin{align*}
\Psi (g(a_1 \cdots a_m), \bm{\tau}) \circ \bm{\sigma} &= \rho_{f, \mu^{[m]} \circ g} \circ (\gamma_{(k_1, p_1)} \otimes \cdots \otimes  \gamma_{(k_m, p_m)})\\
 &\circ (\sigma_1^{\otimes p_1} \otimes \cdots \otimes \sigma_m^{\otimes p_m}) \circ \langle  \Psi (g(a_1 \cdots a_m), \bm{\tau}) \rangle_{k_1^{\times p_1}, \ldots k_m^{\times p_m}} 
\end{align*}
Hence by \eqref{eq sigma determines boldsigma}, we conclude that 

\begin{align*}
(\bm{\tau} \circ \bm{\sigma})_1 \otimes \cdots \otimes (\bm{\tau} \circ \bm{\sigma})_\ell &= \xi_{(g\circ f)(x)} \circ \rho_{f, \mu^{[m]} \circ g} \circ (\gamma_{(k_1, p_1)} \otimes \cdots \otimes  \gamma_{(k_m, p_m)})\\
 &\circ (\sigma_1^{\otimes p_1} \otimes \cdots \otimes \sigma_m^{\otimes p_m}) \circ \langle  \Psi (g(a_1 \cdots a_m), \bm{\tau}) \rangle_{k_1^{\times p_1}, \ldots k_m^{\times p_m}} 
\end{align*}
as we wanted to show.
\end{proof}

\begin{proof}[Proof (of Lemma \ref{lem F(x) bialgebra object})]\label{proof lem F(x) bialg}
As in \ref{ex ComB = fgFMon} all we need to check is that the maps
\begin{align*}
&(\mu, \id): F(x,y) \to F(z) \qquad , \qquad x,y \mapsto z,\\
&(\Delta, \bm{\id}): F(x) \to F(y,z) \qquad , \qquad x \mapsto yz,
\\
&\eta: \bm{1} \to F(x)\\
&\varepsilon :  F(x) \to \bm{1}
\end{align*}
satisfy the bialgebra axioms \eqref{eq algebra axioms} -- \eqref{eq bialgebra axioms} (note that both $\eta$ and $\varepsilon$ lack associated permutations as $\Sigma_0$ is empty). We already checked in \ref{ex ComB = fgFMon} that the equalities hold for the underlying free monoid maps, so we are left to check that the associated permutations coincide. For instance, for the associativity axiom from  Figure \ref{fig bialgebra axioms}(\subref{fig  associativity}) we have that both sides of the equation yield $\id_3 \in \Sigma_3$, since $\langle \id \rangle_{k_1, \ldots}= \id$, $\gamma_{(k,1)}=\id$ and $\xi_{y^2z}= \id$. Perhaps the less obvious verification is the bialgebra axiom from Figure \ref{fig bialgebra axioms}(\subref{fig bialgebra1}). In this case $x,y \mapsto ab$ and the two associated permutations are given by the pair $\bm{\id} = (\id_2, \id_2)$. For both left- and right-hand sides this follows since $(23)(23) = \id_4 = \id_2 \otimes \id_2$. The rest of axioms are also easy to check and are left to the reader.
\end{proof}

We can now wrap up the proof of \ref{thm fundamental theorem}. We are left to check that the ``hat'' construction \eqref{eq def Ahat} defines a two-sided inverse for \eqref{eq bijection alg_B = alg_fgFMonhat}. If $A \in \C$ is a bialgebra, then the composite $$\mathsf{B} \to \widehat{\mathsf{fgFMon}} \overset{\widehat{A}}{\to} \C$$ defines the same bialgebra as $A$. This follows at once from \ref{thm normal form} and the functoriality of $\widehat{A}$. Conversely, if $G: \widehat{\mathsf{fgFMon}} \to \C$ is a functor, let $A \in \C$ be the bialgebra determined by the composite $$\mathsf{B} \to \widehat{\mathsf{fgFMon}} \overset{G}{\to} \C.$$ By construction, the structure maps of $A$ are given by $G(\mu, \bm{\id})$, $G(\Delta, \bm{\id})$, $G(\eta)$ and $G(\varepsilon)$.  Given an arrow $(f, \bm{\sigma})$  in $\widehat{\mathsf{fgFMon}}$, we have that 
\begin{align*}
\widehat{A} (f, \bm{\sigma}) &= G(\mu^{[q_1, \ldots , q_m]}, \bm{\id})  \circ G(P_{\sigma}, \bm{\id}) \circ G(\Delta^{[p_1, \ldots , p_n]}, \bm{\id})\\
&= G((\mu^{[q_1, \ldots , q_m]}, \bm{\id}) \circ (P_{\sigma}, \bm{\id}) \circ (\Delta^{[p_1, \ldots , p_n]}, \bm{\id}))
\end{align*}
so to check that $\widehat{A} =G$ it is enough to show that 
\begin{equation}\label{eq decomposition fsigma}
(f, \bm{\sigma}) = (\mu^{[q_1, \ldots , q_m]}, \bm{\id}) \circ (P_{\sigma}, \bm{\id}) \circ (\Delta^{[p_1, \ldots , p_n]}, \bm{\id})
\end{equation}
where $p_i, q_i$ and $\sigma$ are as described in \eqref{eq def Ahat}. The equality for the underlying monoid map holds since $\Delta^{[p_1, \ldots , p_n]}$ produces as many letters (counting multiplicities) as letters $f(x_1 \cdots x_n)$ has, $P_\sigma$ shuffles them according to the pattern $f(x_1 \cdots x_n)$ and $\bm{\sigma}$ (via $\Psi$) and  $\mu^{[q_1, \ldots , q_m]}$ evaluates the variables producing the pattern. For the equiality at the level of the permutations, by \eqref{eq imposed associativity 1} and \eqref{eq imposed associativity 2} we have that the equality \eqref{eq decomposition fsigma} holds if and only if it holds after precomposition with $(\Delta^{[n]}, \bm{\id})$ and composition with $(\mu^{[m]}, \bm{\id})$. For the left-hand side this yields $$(\mu^{[m]} \circ f \circ  \Delta^{[n]}, \sigma ),$$ where $\sigma = \Psi (f(x_1 \cdots x_n), \bm{\sigma})$ as above. For the right-hand side, we have
\begin{align*}
(\mu^{[m]}, \bm{\id}) &\circ (\mu^{[q_1, \ldots , q_m]}, \bm{\id}) \circ (P_{\sigma}, \bm{\id}) \circ (\Delta^{[p_1, \ldots , p_n]}, \bm{\id}) \circ (\Delta^{[n]}, \bm{\id}) \\
&= (\mu^{[q_1+ \cdots + q_m]}, \bm{\id}) \circ (P_{\sigma}, \bm{\id}) \circ (\Delta^{[p_1 + \cdots + p_n]}, \bm{\id})\\ &= (\mu^{[q_1+ \cdots + q_m]}, \bm{\id}) \circ  (P_{\sigma} \circ \Delta^{[p_1 + \cdots + p_n]}, \bm{\id})\\ 
&= ( \mu^{[q_1+ \cdots + q_m]} P_{\sigma} \circ \Delta^{[p_1 + \cdots + p_n]}, \Psi ((P_{\sigma} \circ \Delta^{[p_1 + \cdots + p_n]})(x), \bm{\id} ))\\
&= ( \mu^{[q_1+ \cdots + q_m]} P_{\sigma} \circ \Delta^{[p_1 + \cdots + p_n]}, \xi^{-1}_{(P_{\sigma} \circ \Delta^{[p_1 + \cdots + p_n]})(x)})\\
&= ( \mu^{[q_1+ \cdots + q_m]} P_{\sigma} \circ \Delta^{[p_1 + \cdots + p_n]}, \xi^{-1}_{x_{\sigma (1)} \cdots x_{\sigma(p_1 + \cdots + p_n)}})\\
&= ( \mu^{[q_1+ \cdots + q_m]} P_{\sigma} \circ \Delta^{[p_1 + \cdots + p_n]}, \sigma)
\end{align*}
as required. This shows that the map \eqref{eq bijection alg_B = alg_fgFMonhat} is indeed a bijection. By \ref{lemma yoneda}, we obtain the equivalence of categories $\widehat{\mathsf{fgFMon}} \simeq \mathsf{B}$. This concludes the proof of \ref{thm fundamental theorem}.

\subsection{A few illustrative examples} To finish off this section, we would like to illustrate the composite law of $\widehat{\mathsf{fgFMon}} $ with a couple of examples.

\begin{example}
Let $(f, \bm{\sigma}): F(x) \to F(a,b)$ and $(g, \tau): F(a,b) \to F(s)$ determined by the following data: $$f(x) = a^2bab \quad, \quad \bm{\sigma} = ((132), (12)) \quad, \quad g(a)=s \quad, \quad g(b)= s^2 \quad, \quad \tau =\id.$$ For any bialgebra $A \in \C$, the map $(f, \bm{\sigma})$  induces the map $$x \mapsto x_{(4)} x_{(1)} x_{(2)} \otimes  x_{(5)} x_{(3)}$$
since $$\Psi (a^2bab, (132), (12)) = \xi_{a^2bab}^{-1} \circ [(132) \otimes (12)] = (34)(132)(45) = (14532).$$ Similarly $(g, \tau)$ induces $a \otimes b \mapsto a b_{(1)} b_{(2)}$. Therefore $\widehat{A}(f, \bm{\sigma}) \circ \widehat{A}(g, \tau)$ must be
$$x \mapsto x_{(4)} x_{(1)} x_{(2)} (x_{(5)} x_{(3)})_{(1)}  (x_{(5)} x_{(3)})_{(2)} =  x_{(5)} x_{(1)} x_{(2)} x_{(6)} x_{(3)}  x_{(7)} x_{(4)}.$$  Let us check that we obtain the same result computing $\widehat{A}( g \circ f, \tau \circ \bm{\sigma})$ applying the composite law \eqref{eq composite law fgFMonhat} (or rather its slightly simplified version \eqref{eq composite law Step 2}). We compute as 
\begin{align*}
\tau \circ \bm{\sigma} &=  \rho_{f,g}  \circ ( \gamma_{(3,1)}\otimes   \gamma_{(2,2)} ) \circ  (\sigma_1  \otimes \sigma_2^{ \otimes 2})  \circ  \langle \id_3 \rangle_{3, 2,2}\\
&= (\langle \xi_{a^2bab} \rangle_{1,1,1,2,2})^{-1} \circ [\id_3 \otimes (23)] \circ [(132)\otimes (12) \otimes (12)]\\
&= (435)(56)(132)(45)(67)\\
&= (1532)(467),
\end{align*}
which yields the same result as above.
\end{example}

\begin{example}
Let $(f, \bm{\sigma}): F(x) \to F(a,b)$, $f(x) = abab$ with $\bm{\sigma} = (\id, (12))$ and $(g, \bm{\tau}): F(a,b) \to F(s,t)$, $g(a)=s, g(b)= ts$ with $\bm{\tau} =((12), \id)$. Let us compute the composite $(g, \bm{\tau}) \circ (f, \bm{\sigma}) $ in $\widehat{\mathsf{fgFMon}}$. To start with, write $$\tau = \Psi (g(ab), (12), \id) = \Psi (sts, (12), \id) = \xi_{sts}^{-1} \circ (12) = (23)(12)= (132).$$ Provided $k_1 = \pi_1 (abab) =2 = \pi_1 (abab) = k_2$ and $p_1 = \pi (s) =1$, $p_2 = \pi (ts)=2$, we have
\begin{align*}
(\bm{\tau} \circ \bm{\sigma})_1 \otimes  (\bm{\tau} \circ \bm{\sigma})_2 &= \xi_{sts^2ts} \circ \rho_{f, \mu  g} \circ (\gamma_{(2, 1)} \otimes  \gamma_{(2,2)}) \circ (\sigma_1 \otimes \sigma_2 \otimes \sigma_2) \circ \langle \tau  \rangle_{2,2,2}\\
&=  (25643) \circ (\langle (23) \rangle_{1,1,2,2})^{-1} \circ (\id_2 \otimes  (23)) \\  &\phantom{--} \circ (\id_2 \otimes (12) \otimes (12)) \circ \langle (132)  \rangle_{2,2,2}\\
&= (25643)(432) (45) (34)(56) (153)(264)\\
&= (25643)(1623) \hspace{5cm} (*)\\
&= (143)(56)\\
&= (143) \otimes (12).
\end{align*}
thence $(\bm{\tau} \circ \bm{\sigma})_1 = (143)$ and $(\bm{\tau} \circ \bm{\sigma})_2 = (12)$.

A bialgebra $A \in \C$ induces the map $\widehat{A}(f, \bm{\sigma})$ mapping $x \mapsto x_{(1)}x_{(3)} \otimes x_{(4)}x_{(2)}$, and similarly $\widehat{A}(g, \bm{\tau})$ maps $a \otimes b  \mapsto b_{(2)}a \otimes b_{(1)}$. Therefore, their composite $\widehat{A}(g, \bm{\tau}) \circ \widehat{A}(f, \bm{\sigma})$ is the map $$x \mapsto (x_{(4)}x_{(2)})_{(2)} x_{(1)}x_{(3)} \otimes (x_{(4)}x_{(2)})_{(1)} = x_{(6)} x_{(3)} x_{(1)} x_{(4)} \otimes x_{(5)} x_{(2)}.$$
We can already see that this agrees with our previous calculation: the rightmost permutation in the equality marked as $(*)$ above is precisely the one that dominates the order of the ``elements'' of the image of $\widehat{A}(g, \bm{\tau}) \circ \widehat{A}(f, \bm{\sigma})$, according to the definition of $\Psi$.
\end{example}

\section{Applications}\label{sec Applications}

We would like to give a few consequences which follow in a rather straightforward way from  theorem \ref{thm fundamental theorem}.

\begin{corollary}
For any symmetric monoidal category $\C$, there is an equivalence of categories $$\mathsf{Alg}_{\mathsf{B}}(\C) \simeq \mathsf{Alg}_{\widehat{\mathsf{fgFMon}}}(\C). $$
\end{corollary}
\begin{proof}
This is immediate from \ref{thm fundamental theorem}.
\end{proof}

\begin{theorem}\label{thm factorisation fgFMonhat}
Every morphism $(f, \bm{\sigma}): F(x_1, \ldots , x_n) \to F(y_1, \ldots , y_m)$ in $\widehat{\mathsf{fgFMon}}$ factors in a unique way as $$f= \mu^{[q_1, \ldots , q_m]} \circ P_{\sigma} \circ \Delta^{[p_1, \ldots , p_n]}  $$ for some (unique) integers $p_1, \ldots , p_n, q_1, \ldots , q_m \geq 0$ such that $\sum_i p_i = \sum_i q_i =:s$, and a (unique) $\sigma \in \Sigma_s$. In particular, these elements are determined as follows: $$ p_i = \pi(f(x_i)) \quad , \quad q_i = \pi_i (f(x_1 \cdots x_n)) \quad , \quad \sigma = \Psi (f(x_1 \cdots x_n), \bm{\sigma}). $$
\end{theorem}
\begin{proof}
This is again immediate from  \ref{thm normal form}, \ref{thm fundamental theorem} and \eqref{eq def Ahat}.
\end{proof}

Now, analogously to \eqref{eq ComA Comb}, there is a canonical faithful functor
$$\mathsf{B} \hooklongrightarrow \mathsf{ComB}  $$
induced by adding the commutativity relation to the symmetric monoidal theory containing the bialgebra structure maps.

\begin{proposition}\label{prop fgFMon hat nohat}
The following diagram is commutative:
$$\begin{tikzcd}
\mathsf{B} \rar{\simeq} \dar &  \widehat{\mathsf{fgFMon}} \dar\\
\mathsf{ComB} \rar{\simeq} & \mathsf{fgFMon}
\end{tikzcd}   $$
where the left vertical functor is the natural inclusion and the right vertical functor forgets the permutations.
\end{proposition}
\begin{proof}
This follows directly from \ref{thm normal form} and \ref{thm fundamental theorem}.
\end{proof}



\begin{proposition}\label{prop fSethat fgFMonhat}
There is a strong monoidal functor $$\widehat{F}: \widehat{\mathsf{fSet}} \to \widehat{\mathsf{fgFMon}} $$ defined as $$\widehat{F} (f, \sigma) := (Ff, \bm{\sigma})$$ where $Ff$ is the monoid map induced by $f$ via the free monoid functor $F$ and $\bm{\sigma}$ is the image of $\sigma$ under the composite $$\Sigma_N \overset{\Psi^{-1}}{\to} \mathcal{W}(k_1, \ldots , k_m) \times \Sigma_{k_1} \times \cdots \times \Sigma_{k_m} \overset{pr}{\to} \Sigma_{k_1} \times \cdots \times \Sigma_{k_m} $$ where $(f, \sigma) : [n] \to [m]$ and $k_i := \# f^{-1}(i)$, which makes the following  diagram commutative,
$$\begin{tikzcd}
\mathsf{A} \dar{\simeq} \rar[hook] &   \mathsf{B} \dar{\simeq}\\
\widehat{\mathsf{fSet}}  \dar \rar{\widehat{F}} & \widehat{\mathsf{fgFMon}}  \dar\\
\mathsf{fSet} \rar{F} & \mathsf{fgFMon}
\end{tikzcd}  $$
where the functors $\widehat{\mathsf{fSet}} \to \mathsf{fSet}$ and  $\widehat{\mathsf{fgFMon}} \to \mathsf{fgFMon}$ are the ones that forget the associated permutations.
\end{proposition}
\begin{proof}
Let us check in first place that the assignment indeed defines a functor. This is certainly true for the underlying monoid maps since $F: \mathsf{fSet} \to \mathsf{fgFMon} $ is already a functor. For $(f, \sigma): [n] \to [m]$ and $(g,  \tau): [m] \to [\ell]$ in $\widehat{\mathsf{fSet}}$,  let us write for convenience $\delta := \sigma \circ \langle \tau \rangle_{k_1, \ldots , k_m}$ for the permutation associated to the composite $(g, \tau)  \circ (f, \sigma)$, and let $\bm{\delta} = (pr \circ \Psi^{-1})(\delta)$. Then
\begin{align*}
(\bm{\tau} \circ \bm{\sigma})_1 \otimes \cdots \otimes (\bm{\tau} \circ \bm{\sigma})_\ell &= \xi_{F(g \circ f)(x_1 \cdots x_n)} \circ( \langle \xi_{(Ff)(x_1 \cdots x_n)} \rangle_{1, \ldots , 1})^{-1} \\
&\phantom{--} \circ (\sigma_1  \otimes  \cdots \otimes \sigma_m) \circ \langle \tau \rangle_{k_1, \ldots , k_m}\\
&= \xi_{x_{gf(1)} \cdots x_{gf(n)}} \circ\xi_{x_{f(1)} \cdots x_{f(n)}}^{-1} \circ (\sigma_1  \otimes  \cdots \otimes \sigma_m)\\ &\phantom{--} \circ \langle \tau \rangle_{k_1, \ldots , k_m}\\
&= \xi_{x_{gf(1)} \cdots x_{gf(n)}} \circ \sigma \circ \langle \tau \rangle_{k_1, \ldots , k_m}\\
&= \xi_{x_{gf(1)} \cdots x_{gf(n)}} \circ \delta\\
&= \delta_1 \otimes \cdots \otimes \delta_\ell,
\end{align*}
where in the third and fifth equality we have used that for a map $(f, \sigma): [n] \to [m]$, the word in $\Psi^{-1}(\sigma)$ can be described as $x_{f(1)} \cdots x_{f(n)}$, as can be easily seen. That $\widehat{F}(\id, \id) = (\id, \bm{\id})$ is straightforward.  The fact that $\widehat{F}$ is strong monoidal follows directly from $F$ being monoidal and the commutative diagram of \ref{lem utilisimo}, which implies that $$(pr \circ \Psi^{-1}) (\sigma \otimes \tau) = ((pr \circ \Psi^{-1})\sigma, (pr \circ \Psi^{-1})\tau),$$ which is precisely the required condition for the associated permutations. Finally, the upper square commutes since the image of $([k] \to [1], \id)$ under $\widehat{F}$ is precisely $(\mu^{[k]}, \id)$ The commutativity of the lower diagram is obvious.
\end{proof}

Putting all the pieces together, we finally have

\begin{theorem}\label{thm main thm 2}
The following cube is commutative:
\begin{equation}
\begin{tikzcd}[row sep=scriptsize,column sep=scriptsize]
& \widehat{\mathsf{fSet}} \arrow[from=dl,"\simeq"] \arrow[rr,"\widehat{F}"] \arrow[dd] & & \widehat{\mathsf{fgFMon}}  \arrow[from=dl,"\simeq"]\arrow[dd] \\
\mathsf{A} \arrow[rr,crossing over,hook]\arrow[dd] & & \mathsf{B} \\
& \mathsf{fSet}\arrow[from=dl, "\simeq"]\arrow[rr, "F \phantom{holaaa}"] &  & \mathsf{fgFMon}\arrow[from=dl,"\simeq"] \\
\mathsf{ComA} \arrow[rr,hook] & & \mathsf{ComB}\arrow[from=uu,crossing over]
\end{tikzcd}
\end{equation}
\end{theorem}
\begin{proof}
This follows at once from \ref{lem fSet fgFMon}, \ref{thm fundamental theorem}, \ref{prop fgFMon hat nohat} and \ref{prop fSethat fgFMonhat}. 
\end{proof}

If $\mathsf{C}$  (resp. $\mathsf{CocomC}$) is the PROP for (cocommutative) coalgebras, by reversing arrows we directly obtain

\begin{corollary}\label{cor main thm for coalg}
There are monoidal equivalences $$\mathsf{C} \overset{\simeq}{\to}\widehat{\mathsf{fSet}}^{op} \qquad , \qquad \mathsf{B} \overset{\simeq}{\to} \widehat{\mathsf{fgFMon}}^{op} $$
making the following cube commutative,
\begin{equation}
\begin{tikzcd}[row sep=scriptsize,column sep=scriptsize]
& \widehat{\mathsf{fSet}}^{op} \arrow[from=dl,"\simeq"] \arrow[rr,"\widehat{F}^{op}"] \arrow[dd] & & \widehat{\mathsf{fgFMon}}^{op}  \arrow[from=dl,"\simeq"]\arrow[dd] \\
\mathsf{C} \arrow[rr,crossing over,hook]\arrow[dd] & & \mathsf{B} \\
& \mathsf{fSet}^{op}\arrow[from=dl, "\simeq"]\arrow[rr, "F^{op} \phantom{holaaa}"] &  & \mathsf{fgFMon}^{op}\arrow[from=dl,"\simeq"] \\
\mathsf{CocomC} \arrow[rr,hook] & & \mathsf{CocomB}\arrow[from=uu,crossing over]
\end{tikzcd}
\end{equation}
\end{corollary}

\subsection{Related work}
Our work could be seen as the bialgebra case of what is called in the literature as the \textit{acyclic Hopf word problem} \cite{thomas}. It gives a complete answer to the question of whether two different combinations of the bialgebra structure maps represent the same morphism in $\mathsf{B}$. More precisely, it asks whether given two arrows $f_1,f_2: n \to m$ in the PROP $\mathcal{B}$ as defined just before \ref{lem elements in normal form}, it is true that $\pi (f_1)= \pi (f_2)$, where $\pi: \mathcal{B} \to \mathsf{B}$ is the canonical projection. Indeed if $(f_1', \bm{\sigma})$, $(f_2', \bm{\tau})$ are the images of $\pi (f_1), \pi (f_2)$ under the equivalence $\mathsf{B} \overset{\simeq}{\to}  \widehat{\mathsf{fgFMon}} $, then \ref{thm factorisation fgFMonhat} ensures that $\pi (f_1)= \pi (f_2)$ if and only if the following equalities hold:
$$  \pi(f_1'(x_i)) =  \pi(f_2'(x_i)) \quad , \quad  \pi_j (w_1) = \pi_j(w_2) \quad , \quad \Psi (w_1, \bm{\sigma})= \Psi (w_2, \bm{\tau}) $$ for all $1 \leq i \leq n$, $1\leq j \leq m$, where $w_1 := f_1' (x_1 \cdots x_n)$ and $w_2 := f_2' (x_1 \cdots x_n)$.

In \cite{thomas} (and in \cite{kuperberg} in a different context) this problem is used to produce an invariant of decorated 3-manifolds.


\subsection{Future work} The commutative cube from \ref{thm main thm} is only a part of the larger solid diagram shown below.

\begin{equation}
\begin{tikzcd}[row sep=scriptsize,column sep=scriptsize]
& \widehat{\mathsf{fSet}} \arrow[from=dl,"\simeq"] \arrow[rr,"\widehat{F}"] \arrow[dd] & & \widehat{\mathsf{fgFMon}}  \arrow[from=dl,"\simeq"]\arrow[dd] \arrow[rr, dashed] & & \widehat{\mathsf{fgFGp}}  \arrow[from=dl,"\simeq", dashed]\arrow[dd, dashed] \\
\mathsf{A} \arrow[rr,crossing over,hook]\arrow[dd] & & \mathsf{B}  \arrow[rr,crossing over,hook] & & \mathsf{H} \\
& \mathsf{fSet}\arrow[from=dl, "\simeq"]\arrow[rr, "F \phantom{holaaa}"]  &  & \mathsf{fgFMon}\arrow[from=dl,"\simeq"] \arrow[rr, "G \phantom{holaaa}"] \arrow[dd] &  & \mathsf{fgFGp}\arrow[from=dl,"\simeq"] \arrow[dd] \\
\mathsf{ComA} \arrow[rr,hook] & & \mathsf{ComB}\arrow[from=uu,crossing over]  \arrow[rr,crossing over,hook] & & \mathsf{ComH}\arrow[from=uu,crossing over]\\ 
&  &  & \mathsf{fgFCMon}\arrow[from=dl,"\simeq"] \arrow[rr, "K \phantom{holaaa}"] &  & \mathsf{fgFAb}\arrow[from=dl,"\simeq"] \\
& & \mathsf{ComCocomB} \arrow[from=uu,crossing over]  \arrow[rr,hook] & & \mathsf{ComCocomH} \arrow[from=uu,crossing over]
\end{tikzcd}
\end{equation}


Let us say a few words about the this diagram. Here $\mathsf{H}$ (resp. $\mathsf{ComH}$) stands for the PROP for (commutative) Hopf algebras. All functors in the diagram are bijective-on-objects. Horizontal functors are faithful,  vertical functors are full, and perpendicular (with respect to the page) functors are equivalences of categories. The functor $K$ is the restriction of the \textit{Grothendieck construction} $K: \mathsf{ CMon} \to \mathsf{Ab} $, that is, the left-adjoint to the forgetful $\mathsf{Ab} \to \mathsf{ CMon}$, that topologists should be familiar with. The functor  $G$ is its non-commutative counterpart, the restriction of  the left-adjoint to the forgetful $\mathsf{Gp} \to \mathsf{Mon}$, which is known as the \textit{group completion} or \textit{universal enveloping group} of a monoid, $G: \mathsf{Mon} \to \mathsf{Gp}$. The vertical functors $\mathsf{Mon} \to \mathsf{CMon}$ and $\mathsf{Gp} \to \mathsf{Ab}$ are the abelianisation functors.

In this situation, we conjecture the existence of a category $\widehat{\mathsf{fgFGp}}$ built out of $\mathsf{fgFGp}$ in an analogous way to our construction of $\widehat{\mathsf{fgFMon}}$, together with a monoidal equivalence $\mathsf{H} \simeq \widehat{\mathsf{fgFGp}}$, a functor $\widehat{G}: \widehat{\mathsf{fgFMon}} \to \widehat{\mathsf{fgFGp}}$ lifting $G$ and a full (forgetful) functor $\widehat{\mathsf{fgFGp}} \to \mathsf{fgFGp}$. This will be the subject of a forthcoming paper.

\bibliographystyle{alpha}
\bibliography{bibliografia}

\end{document}